\documentclass[aos,MSNbibl,nameyear,dvips]{arximspdf}
\usepackage{mathbh}

\doi{10.1214/11-AOS912}
\volume{39}
\issue{5}
\pubyear{2011}
\firstpage{2557}
\lastpage{2584}

\makeatletter

\newcommand{\N}{\mathbb{N}}
\newcommand{\R}{\mathbb{R}}
\newcommand{\LL}{\mathbb{L}}

\newcommand{\E}{E}
\newcommand{\Proba}{P}
\newcommand{\tr}{\operatorname{tr}}
\newcommand{\ud}{d}
\newcommand{\phispan}{{\langle\phi\rangle}}
\newcommand{\covmat}{\Sigma}

\newtheorem{theorem}{Theorem}
\newtheorem{lem}{Lemma}
\newtheorem{corol}{Corollary}
\newtheorem{prop}{Proposition}
\newproclaim{df}{Definition}

\makeatother

\begin{document}
\begin{frontmatter}

\title{Bernstein--von Mises theorems for Gaussian regression with
increasing number of regressors}
\runtitle{Bernstein--von Mises theorems for Gaussian regression}

\begin{aug}
\author[A]{\fnms{Dominique} \snm{Bontemps}\corref{}\ead[label=e1]{dominique.bontemps@gmail.com}}
\runauthor{D. Bontemps}
\affiliation{Universit\'{e} Paris-Sud}
\address[A]{Laboratoire de Math\'{e}matiques d'Orsay\\
Universit\'{e} Paris-Sud\\
UMR8628, B\^{a}t. 425\\
F-91405, Orsay\\
France\\
\printead{e1}} 
\end{aug}

\received{\smonth{9} \syear{2010}}
\revised{\smonth{6} \syear{2011}}

%
\begin{abstract}
This paper brings a contribution to the Bayesian theory of
nonparametric and semiparametric estimation. We are interested in the
asymptotic normality of the posterior distribution in Gaussian linear
regression models when the number of regressors increases with the
sample size. Two kinds of Bernstein--von Mises theorems are obtained in
this framework: nonparametric theorems for the parameter itself, and
semiparametric theorems for functionals of the parameter. We apply them
to the Gaussian sequence model and to the regression of functions in
Sobolev and $C^{\alpha}$ classes, in which we get the minimax
convergence rates. Adaptivity is reached for the Bayesian estimators of
functionals in our applications.
\end{abstract}

%

\end{frontmatter}

\section{Introduction} \label{secIntro}

To estimate a parameter of interest in a statistical mo\-del, a~Bayesian
puts a prior distribution on it and looks at the posterior
distribution, given the observations. A Bernstein--von Mises theorem is
a~result giving conditions under which the posterior distribution is
asymptotically normal, centered at the maximum likelihood estimator
(MLE) of the model used, with a variance equal to the asymptotic
frequentist variance of the MLE. 
Other centering can be used; see, for instance, \citet
{VanderVaartAsympStats98}, page 144, after the proof of Lemma 10.3.

Such an asymptotic posterior normality is important because it allows
the construction of approximate credible regions, based on the
posterior distribution, which retain good frequentist properties. In
particular, the Monte Carlo Markov chain algorithms (MCMC) make
feasible the construction of Bayesian confidence regions in complex
models, for which frequentist confidence regions are difficult to
build; however, Bernstein--von Mises theorems are difficult to derive in
complex models.

Note that the Bernstein--von Mises theorem also has links with
information theory [see \citet{ClarkeB90} and \citet
{ClarkeGhosal2010}].

For parametric models, the Bernstein--von Mises theorem is a well-known
result, for which we refer to \citet{VanderVaartAsympStats98}. In
nonparametric models (where the parameter space is infinite-dimensional
or growing) and semiparametric models (when the parameter of interest
is a finite-dimensional functional of the complete infinite-dimensional
parameter), there are still relatively few asymptotic normality
results. \citet{Freedman99Wald} gives negative results, and we recall
some positive ones below. However, many recent papers deal with the
convergence rate of posterior distributions in various settings, which
is linked with the model complexity: we refer to \citet
{GhosalGhoshVanderVaart00}, \citet{ShenWasserman2001} as early representatives
of this school.

Nonparametric Bernstein--von Mises theorems have been developed for
models based on a sieve approximation, where the dimension of the
parameter grows with the sample size. In particular, two situations
have been studied: regression models in \citet{Ghosal99}; exponential
models in \citet{Ghosal00}, \citet{ClarkeGhosal2010} and
\citet
{BoucheronGassiat09} (this last one deals with the discrete case, when
the observations follow some unknown infinite multinomial
distribution).

In semiparametric frameworks the asymptotic normality has been obtained
in several situations. \citet{KimLee2004BvM} and \citet
{Kim2006BvM} study
the nonparametric right-censoring model and the proportional hazard
model. \citet{Castillo08BvM} obtains Bernstein--von Mises theorems for
Gaussian process priors, in the semiparametric framework where the
unknown quantity is $(\theta, f)$, with $\theta$ the parameter of
interest and $f$ an infinite-dimensional nuisance parameter. See also
\citet{Shen2002}. \citet{RivoirardRousseau2009} obtain the
Bernstein--von
Mises theorem for linear functionals of the density of the
observations, in the context of a sieve approximation: sequences of
spaces with an increasing dimension $k_n$ are used to approximate an
infinite-dimensional space. These authors achieve also the frequentist
minimax estimation rate for densities in specific regularity classes
with a deterministic (nonadaptive) value of the dimension $k_n$.

Here we obtain nonparametric and semiparametric Bernstein--von Mises
theorems in a Gaussian regression framework with an increasing number
of regressors. We address two challenging problems. First, we try to
understand better when the Bernstein--von Mises theorem holds and when
it does not. In the latter case the Bayesian credible sets no longer
preserve their frequentist asymptotic properties. Second, we look for
adaptive Bayesian estimators in our semiparametric settings.

Our nonparametric results cover the case of a specific Gaussian prior,
and the case of more generic smooth priors. They are said to be
nonparametric because we use sieve priors, that is, the dimension of
the parameter grows. These results improve on the preceding ones by
\citet{Ghosal99} which did not suppose the normality of the errors but
imposed other conditions, in particular, on the growth rate of the
number of regressors. We apply our results to the Gaussian sequence
model, as well as to periodic Sobolev classes and to regularity classes
$C^{\alpha}[0, 1]$ in the context of the regression model (using,
resp., trigonometric polynomials and splines as regressors). In
all these situations we get the asymptotic normality of the posterior
in addition to the minimax convergence rates, with appropriate
(nonadaptive) choices of the prior. We also show that for some priors
known to reach this convergence rate, the Bernstein--von Mises theorem
does not hold.

We derive also semiparametric Bernstein--von Mises theorems for linear
and nonlinear functionals of the parameter. The linear case is an
immediate corollary of the nonparametric theorems and does not need any
additional conditions. We apply these results to the periodic Sobolev
classes to estimate a linear functional and the $L^2$ norm of the
regression function $f$ when it is smooth enough, and in both cases we
are able to build an adaptive Bayesian estimator which achieves the
minimax convergence rate 
in all classes of the collection,
in addition to the asymptotic normality.

The paper is organized as follows. We present the framework in
Section~\ref{secFramework}. Section \ref{secnonparametric} states the
nonparametric Bernstein--von Mises theorems, for Gaussian and
non-Gaussian priors. In Section \ref{secsemiparametric} we derive the
semiparametric Bernstein--von Mises theorems for linear and nonlinear
functionals of the parameter. Then in Section \ref{secapplications} we
give applications to the Gaussian sequence model, and to the regression
of a function in a Sobolev and $C^{\alpha}[0, 1]$ class. In Section
\ref
{secproofs} the nonparametric and semiparametric Bernstein--von Mises
theorems are proved. The appendices contain various technical tools
used in the main analysis; the appendices can be found in the
supplemental article [\citet{BvMGRegSupplement}].

\section{Framework} \label{secFramework}

We consider a Gaussian linear regression framework. For any $n\geq1$,
our observation $Y=(Y_1, \ldots, Y_n)\in\R^n$ is a Gaussian random vector
%
%
\begin{equation} \label{eqvector-model}
Y = F + \varepsilon,
\end{equation}
where the vector of errors $\varepsilon=(\varepsilon_1, \ldots,
\varepsilon_n) \sim\mathcal{N}(0, \sigma_n^2 I_n)$, with $I_n$ the
\mbox{$n\times n$} identity matrix, and the mean vector $F$ belongs to $\R^n$.
Note that the dimension of $Y$ is the sample size $n$, and that $\sigma
_n^2$ is known but may depend on $n$. Let $F_0$ be the true mean vector
of $Y$ with distribution $\mathcal{N}(F_0, \sigma_n^2 I_n)$.
Probability expectations under $F_0$ are denoted $\Proba_{F_0}$ and
$\E$.


Let $\phi_1, \ldots, \phi_{k_n}$ a collection of $k_n$ linearly
independent regressors in $\R^n$, where $k_n\leq n$ grows with $n$. We
gather these regressors in the $n\times k_n$-matrix $\Phi$ of rank
$k_n$, and $\phispan= \{ \Phi\theta\dvtx\theta=(\theta_1, \ldots,
\theta_{k_n}) \in\R^{k_n}\}$ denotes their linear span. The
Bernstein--von Mises theorems will be stated in association with~%
$\phispan$, the vector space of possible mean vectors in the model,
which is possibly misspecified. We denote by $\Proba_\theta$ the
probability distribution of a random variable following $\mathcal
{N}(\Phi\theta, \sigma_n^2 I_n)$ and $\E_\theta$ the associated
expectation.

As examples, we present three different settings, each with its own
collection of regressors. In Section \ref{secapplications} the
Bernstein--von Mises theorems are applied to each of these frameworks:
\begin{longlist}[(1)]
\item[(1)] \textit{The Gaussian sequence model.}
Our first application concerns the Gaussian sequence model, which is
also equivalent to the white noise model [see \citet{Massart2007},
Chapter 4, e.g.]. We consider the infinite-dimensional setting
%
%
\begin{equation} \label{eqwhite-noise}
Y_j = \theta_j^0 + \frac{1}{\sqrt{n}} \xi_j,\qquad j\geq1,
\end{equation}
where the random variables $\xi_j, j\geq1$ are independent and have
distribution~$\mathcal{N}(0,1)$. Projecting on the first $k_n$
coordinates with $k_n \leq n$, we retrieve our model (\ref
{eqvector-model}) with 
$\theta_0 = (\theta_j^0)_{1\leq j\leq k_n}$, $\sigma_n = 1/\sqrt{n}$
and $\Phi^T \Phi= I_{k_n}$.

\item[(2)] \textit{Regression of a function in a Sobolev class.}
Let $f\dvtx[0, 1] \rightarrow\R$ be a~function in $\LL^2([0,1])$. We
observe realizations of random variables
%
%
\begin{equation} \label{eqregression-0-1}
Y_i = f(i/n) + \varepsilon_i
\end{equation}
for $1 \leq i \leq n$, where the errors $\varepsilon_i$ are i.i.d.
$\mathcal{N}(0, \sigma_n^2)$ and $\sigma_n$ does not depend on $n$.

We denote by $(\varphi_j)_{j\geq1}$ the Fourier basis
%
%
\begin{eqnarray} \label{eqFourier-basis}
\varphi_1 &\equiv& 1, \nonumber\\
\varphi_{2 m}(x) &=& \sqrt{2} \cos(2 \pi m x) \qquad\forall m\geq1, \\
\varphi_{2 m+1}(x) &=& \sqrt{2} \sin(2 \pi m x) \qquad\forall m\geq
1.\nonumber
\end{eqnarray}

In conjunction with the regular design $x_i = i/n$ for $1\leq i\leq
n$, this gives the collection of regressors
\[
\phi_j = \bigl(\varphi_j(i/n)\bigr)_{1\leq i\leq n},\qquad 1\leq j\leq
k_n.
\]

In practice, we suppose that $f$ belongs to one of the periodic
Sobolev classes:
\begin{df} \label{dfSobolev}
Let $\alpha>0$ and $L>0$. Let $(\varphi_j)_{j\geq1}$ denote the
Fourier basis~(\ref{eqFourier-basis}). We define the Sobolev class
$\mathcal{W}(\alpha, L)$ as the collection of all functions $f = \sum
_{j=1}^{\infty} \theta_j \varphi_j$ in $\LL^2([0, 1])$ such that
$\theta
=(\theta_j)_{j\geq1}$ is an element of the ellipsoid of $\ell^2(\N)$,
\[
\Theta(\alpha, L) = \Biggl\{ \theta\in\ell^2(\N)\dvtx\sum_{j=1}^{\infty}
a_j^2 \theta_j^2 \leq\frac{L^2}{\pi^{2\alpha}} \Biggr\},
\]
where
%
%
\begin{equation}
a_j = \cases{
j^\alpha, &\quad if $j$ is even;\cr
(j-1)^\alpha, &\quad if $j$ is odd.}
\end{equation}
\end{df}

\item[(3)] \textit{Regression of a function in $C^\alpha[0,1]$.}
Let $\alpha>0$, and $f\in C^\alpha[0,1]$. This means that $f$ is
$\alpha_0$ times continuously differentiable with $\|f\|_{\alpha}
<\infty$, $\alpha_0$ being the greatest integer less than $\alpha$ and
the seminorm \mbox{$\|\cdot\|_{\alpha}$} being defined by
\[
\|f\|_{\alpha} = \sup_{x\neq x'} \frac{|f^{(\alpha
_0)}(x)-f^{(\alpha_0)}(x')|}{|x-x'|^{\alpha-\alpha_0}}.
\]
Consider a design $(x_i^{(n)})_{n\geq1, 1\leq i \leq n}$, not
necessarily uniform. Here $F_0$ is the vector $(f(x_i^{(n)}
))_{1\leq i \leq n}$. Once again we suppose that $\sigma_n=\sigma$
does not depend on $n$.

Fix an integer $q\geq\alpha$, and let $K=k_n+1-q$. Partition the
interval $(0, 1]$ into $K$ subintervals $( (j-1)/K, j/K]$
for $1\leq j \leq K$. We want to perform the regression of $f$ in the
space of splines of order $q$ defined on that partition, and use the
$B$-splines basis $(B_j)_{1\leq j\leq k_n}$ [see, e.g.,
\citet
{deBoor78}]. Our collection of regressors is $\phi_j = (B_j
(x_i^{(n)}))_{1\leq i\leq n}$, for $1\leq j\leq k_n$.
\end{longlist}

For any value of $n\geq1$, let $\widetilde{W}$ be a prior distribution
on $\R^{k_n}$ and, for $F=\Phi\theta$, let $W$ be the prior
distribution on $F\in\R^n$ obtained from $\widetilde{W}$ on~$\theta$.
Its support is included in $\phispan$. 
Let $P^W$ denote the marginal distribution of $Y$ under prior $W$, and
$W(\ud G(F) | Y)$ denote the posterior distribution of a functional
$G(F)$. Note that everything depends on $n$ ($W$, e.g., is a~%
distribution on $\R^n$) even if we do not use $n$ as an index to
simplify our notation.

Both the parametrization by $\theta$ and the corresponding collection
of regressors $\phi_1, \ldots, \phi_{k_n}$ are arbitrary: what matters
is the posterior distribution of~$F$ and this depends on the space
$\phispan$, not on the basis used to parametrize it. The span $\phispan
$ is characterized by the matrix $\covmat= \Phi(\Phi^T \Phi)^{-1} \Phi
^T$ of the orthogonal projection onto~$\phispan$. 

The prior $W$ is a sieve prior: that is, its support comes from a
finite-dimensional model whose dimension $k_n$ grows with $n$. The
collection of growing models $\phispan$ (the sieve) can be seen as an
approximation framework, each model being possibly misspecified. There
is no true parameter in our setting: the true mean vector $F_0$ may
fall outside $\phispan$ and correspond to none of the possible values
of $\theta$. There is then a bias which has to be dealt with, linked to
the choice of the cutoff $k_n$.

When dealing with Bernstein--von Mises results, the question of the
asymptotic centering point arises. In nonparametric models constructed
on an infinite-dimensional parameter, there is no definition of a MLE;
what the natural centering for a Bernstein--von Mises theorem should be
in such situations is not clear. In the model $\phispan$, the
orthogonal\vspace*{1pt} projection $Y_{\phispan}=\covmat Y$ of $Y$ is also the MLE
of $F_0$. We set $\theta_Y=(\Phi^T \Phi)^{-1} \Phi^T Y$ its associated
parameter. Let also $F_{\phispan} = \Phi\theta_0$ be the projection of
$F_0$ on $\phispan$, with $\theta_0=(\Phi^T \Phi)^{-1} \Phi^T F_0$.
Now, $F_0-F_\phispan$ corresponds to the bias introduced by the use of
the model~$\phispan$, and $F_\phispan$ is the centering point of the
distribution of the MLE $Y_{\phispan}$ under~$\Proba_{F_0}$:
\[
Y_{\phispan} \sim\mathcal{N}\bigl(F_{\phispan}, \sigma_n^2 \covmat
\bigr).
\]
Although the MLE is naturally defined \textit{in the sieve $\phispan$},
it heavily depends on the choice of $\phispan$. Therefore, the
Bernstein--von Mises theorems we establish depend on the choice of the
sieve the prior distribution is built on. 

\section{Nonparametric Bernstein--von Mises theorems} \label{secnonparametric}

The proofs of our nonparametric results are delayed to Section \ref
{secproofs}.

\subsection{With Gaussian priors} \label{secallnormal}

We consider here a centered, normal prior distribution $W$ which is
isotropic on $\phispan$, so that $W=\mathcal{N}(0, \tau_n^2 \covmat
)$ for some sequence~$\tau_n$. $\tau_n$ is a scale parameter, and
essentially the only assumption needed in this case is that $\tau_n$ is
large enough 
as $n$ grows. Let $\| Q - Q' \|_{\mathrm{TV}}$ denote the total
variation norm between two probability distributions $Q$ and~$Q'$.
\begin{theorem} \label{thmall-Gaussian-BvM}
Assume that $\sigma_n=o(\tau_n)$, $\|F_0\| = o(\tau_n^2/\sigma_n)$ and
$k_n =\break o(\tau_n^4/\sigma_n^4)$.
Then
\[
\E\bigl\| W(\ud F | Y) - \mathcal{N}\bigl(Y_{\phispan}, \sigma_n^2
\covmat\bigr) \bigr\|_{\mathrm{TV}} \rightarrow0 \qquad\mbox{as }
n\rightarrow\infty.
\]
\end{theorem}

In terms of $\theta$ instead of $F$, an equivalent statement is
\[
\E\bigl\| \widetilde{W}(\ud\theta| Y) - \mathcal{N}(\theta_Y,
\sigma_n^2 (\Phi^T \Phi)^{-1} ) \bigr\|_{\mathrm{TV}}
\rightarrow0 \qquad\mbox{as } n\rightarrow\infty.
\]

Theorem \ref{thmall-Gaussian-BvM} does not deal with the modeling bias
introduced by taking a prior restricted to $\phispan$. This is an
important question in nonparametric statistics, and $k_n$ has to be
chosen in order to achieve a satisfactory bias-variance trade-off. 

As an example, let us consider a typical regression framework with
$F_0=(f_0(x_i))_{1\leq i\leq n}$, where $f_0$ is some function and
$(x_i)_{1\leq i\leq n}$ some design. If $\sigma_n$ does not depend on
$n$, both conditions $\|F_0\| = o(\tau_n^2/\sigma_n)$ and $k_n =
o(\tau_n^4/\sigma_n^4)$ are satisfied if $f_0$ is bounded and
$n^{1/4}=o(\tau_n)$. These\vspace*{1pt} conditions can be read in
another way: $\tau _n^4$ must be large enough with respect to $\|F_0\|$
and $k_n$.

\subsection{With smooth priors} \label{secsmooth}

We consider now more general priors. To understand better the
conditions we use, we need to look at the mechanics of the
Bernstein--von Mises theorem.\vadjust{\goodbreak}

Behind a Bernstein--von Mises theorem there is a LAN structure: the
log-likelihood admits a quadratic expansion near the MLE. Since the
posterior density is proportional to the product of the prior density
and the likelihood, the prior has to be locally constant to let the
likelihood alone influence the posterior and produce the Gaussian
shape. To prove a Bernstein--von Mises theorem, we look for a subset
which is simultaneously (1) large enough, so that the posterior will
concentrate on it, and (2) small enough, so that we can find
approximately constant priors on it. The larger the dimension of the
model is, the more difficult it is to combine these two requirements,
and the more difficult it is to obtain a Bernstein--von Mises theorem.

The geometry of the subsets are naturally suggested by the normal
distribution we are looking for. For $M>0$, consider the ellipsoid
%
%
\begin{equation} \label{eqdef-ellipsoid}
\mathcal{E}_{\theta_0, \Phi}(M) = \{ \theta\in\R^{k_n}\dvtx(\theta
-\theta_0)^T \Phi^T\Phi(\theta-\theta_0) \leq\sigma_n^2 M \}.
\end{equation}
\begin{theorem} \label{thmsmooth-gaussian-BvM}
Suppose that $W$ is induced by a distribution $\widetilde{W}$ on
$\theta$ admitting a density $w(\theta)$ with respect to the Lebesgue
measure. If there exists a sequence $(M_n)_{n\geq1}$ such that:
\renewcommand\thelonglist{(\arabic{longlist})}
\renewcommand\labellonglist{\thelonglist}
\begin{longlist}
\item\label{hypcontinuous-w} $ \sup_{\|\Phi h\|^2 \leq
\sigma_n^2 M_n, \|\Phi g\|^2 \leq\sigma_n^2 M_n} \frac{w(\theta
_0+h)}{w(\theta_0+g)} \rightarrow1$
as $n\rightarrow\infty$,
\item\label{hypk-n-M-n} $k_n \ln k_n = o(M_n)$,
\item\label{hypweight} $ \max(0, \ln(\frac
{\sqrt{\det(\Phi^T\Phi)}}{\sigma_n^{k_n} w(\theta_0)} )) = o(M_n)$,
\end{longlist}
then
\[
\E\bigl\| W(\ud F | Y) - \mathcal{N}\bigl(Y_{\phispan}, \sigma_n^2
\covmat\bigr) \bigr\|_{\mathrm{TV}} \rightarrow0 \qquad\mbox{as }
n\rightarrow\infty.
\]
\end{theorem}

With condition \ref{hypcontinuous-w} below we ask for a sufficiently
flat prior $\widetilde{W}$ in an ellipsoid $\mathcal{E}_{\theta_0,
\Phi
}(M_n)$. Condition \ref{hypk-n-M-n} ensures, in particular, that the
weight the normal distribution puts on $\mathcal{E}_{\theta_0, \Phi
}(M_n)$ in the limit goes to $1$. Condition \ref{hypweight} makes
quantities linked to the volume of $\mathcal{E}_{\theta_0, \Phi}(M_n)$
appear and guarantees that it has enough prior weight. This kind of
assumption is common in the literature dealing with the concentration
of posterior distributions; see, for instance, \citet
{GhosalGhoshVanderVaart00}.


Several of our applications illustrate that priors known to induce the
posterior minimax convergence rate may not be flat enough to get the
Gaussian shape with the asymptotic variance $\sigma_n^2 \covmat$.

An important remark is the following: condition \ref{hypk-n-M-n} does
not really limit the growth rate of $k_n$. Read in conjunction with the
other two conditions, we see that a flatter prior distribution will
permit us to take $M_n$ larger. Thus, the only condition on the growth
rate of $k_n$ is $k_n\leq n$.

Note that Theorem \ref{thmsmooth-gaussian-BvM} is not a generalization
of Theorem \ref{thmall-Gaussian-BvM}: Theorem \ref
{thmall-Gaussian-BvM} is more powerful for isotropic Gaussian priors.
Consider again the regression framework with $F_0=(f_0(x_i)
)_{1\leq i\leq n}$, where $f_0$ is a bounded function and $(x_i)_{1\leq
i\leq n}$ is some design. Suppose\vspace*{1pt} that $\sigma_n$ does not depend on
$n$, and take $k_n=n$ and $W=\mathcal{N}(0, \tau_n^2 \covmat
)$. Then the conditions of Theorem \ref{thmall-Gaussian-BvM} are
satisfied as soon as $n^{1/4}=o(\tau_n)$, but with Theorem \ref
{thmsmooth-gaussian-BvM} we need $n \ln n = o(\tau_n^2)$.

Our main applications, to the Gaussian sequence model and to the
regression model using trigonometric polynomials and splines, are
developed in Section \ref{secapplications}. We now present two remarks
about the parametric case and the comparison with the pioneer work of
\citet{Ghosal99}.

\subsubsection*{The parametric case} Consider the regression of a function
$f$ defined on $[0, 1]$, with a fixed number $k$ of regressors. Set a
design $(x_i^{(n)})_{n\geq1, 1\leq i\leq n}$, with $x_i^{(n)}\in
[(i-1)/n, i/n]$ for any $n \geq1$, and $F_0=(f(x_i^{(n)})
)_{1\leq i\leq n}$. Choose a finite number of piecewise continuous and
linearly independent regressors $(\varphi_j)_{1\leq j\leq k}$ on $[0,
1]$, and set $\phi_j = (\varphi_j(x_i^{(n)})
)_{1\leq i \leq n}$ for $1\leq j\leq k$. Assume that~$f$, $k_n=k$,
$\sigma_n=\sigma$ and $\widetilde{W}$ do not depend on $n$.

We would like to compare Theorem \ref{thmsmooth-gaussian-BvM} with the
usual Bernstein--von Mises theorem for parametric models applied to such
a regression framework. In that setting, let us suppose that $w$ is
continuous and positive, and that~$f$ is bounded. Then condition \ref
{hypcontinuous-w} becomes $M_n = o(n)$, while condition \ref
{hypweight} reduces to $\ln n = o(M_n)$. Clearly, there exist such
sequences $(M_n)_{n\geq1}$, so Theorem \ref{thmsmooth-gaussian-BvM}
applies. Here the rescaling by $\sqrt{n}$ of the Bernstein--von Mises
theorem for parametric models is hidden in the asymptotic posterior
variance $\sigma^2 (\Phi^T\Phi)^{-1}$ of the parameter $\theta$.
Indeed, $(1/n)$ $\Phi^T\Phi$ is a Riemann sum and converges toward the
Gramian matrix of the collection $(\varphi_j)_{1\leq j\leq k}$ in~$\LL
^2([0,1])$.
\begin{pf}
$\!\!$We have\vspace*{1pt} $\|\Phi\theta_0\| \leq\|F_0\| \leq\sqrt{n} \|f\|_\infty$,
and $\|\theta_0\|^2 \leq\| (\Phi^T\Phi)^{-1}\| \cdot\|\Phi
\theta_0\|^2 \leq\| n (\Phi^T\Phi)^{-1}\| \|f\|_\infty^2$.
$(1/n)$ $\Phi^T\Phi$ converges toward the Gramian matrix of the
collection $(\varphi_j)_{1\leq j\leq k}$ in $\LL^2([0,1])$, and its
smallest eigenvalue is lower bounded for $n$ large enough. Therefore,
$\theta_0$ is bounded, and we can consider it lies in some compact set
on which $w$ is uniformly continuous and lower bounded by a positive
constant. The rest follows.
\end{pf}

\subsubsection*{Comparison with Ghosal's conditions} 
The Bernstein--von Mises theorem in a regression setting when the number
of parameters goes to infinity has been first studied by 
\citet{Ghosal99} as an early step in the development of frequentist
nonparametric Bayesian theory. In his paper the errors $\varepsilon_i$
are not supposed to be Gaussian. Under the Gaussianity assumption we
get improved results, which means that we have a nontrivial
generalization of the \citet{Ghosal99} conditions in the case of
Gaussian errors. In particular, our condition for the prior smoothness
is simpler, and the growth rate of the dimension $k_n$ is much less constrained:
\begin{itemize}
\item\citet{Ghosal99} does not admit a modeling bias between $F_0$ and
$\Phi\theta_0$. In the present work the normality of the errors
permits us to take $F_0\neq\Phi\theta_0$ without any cost, as it
appears in the core of the proof (Lemma~\ref
{leminvariant-translation}). The possibility of considering
misspecified models is an important improvement.
\item In \citet{Ghosal99} $\sigma_n$ is constant, which does not allow
the application to the Gaussian sequence model.
\item\citet{Ghosal99} restricts the growth of the dimension $k_n$ to
$k_n^4\ln k_n = o(n)$ (see below). In our setting we only require
$k_n\leq n$. With Ghosal's condition we could not have obtained the
applications to the Gaussian sequence model or to the regression model
for Sobolev or $C^{\alpha}$ classes.
\end{itemize}

Let $\delta_n^2=\|(\Phi^T\Phi)^{-1}\|$ be the operator norm of
$(\Phi
^T\Phi)^{-1}$ for the $\ell^2$ metric, and let $\eta_n^2$ be the
maximal value on the diagonal of $\covmat$. With our notation, the last
two assumptions of \citet{Ghosal99} become:
\begin{longlist}[(A3)]
\item[(A3)] 
There exists $\eta_0>0$ such that $w(\theta_0)>\eta_0^{k_n}$. Moreover,
%
%
\begin{equation} \label{eqhyp-24}
|{\ln w}(\theta) - \ln w(\theta_0)| \leq L_n(C) \|\theta-\theta_0\|,
\end{equation}
whenever $\|\theta-\theta_0\| \leq C \delta_n k_n \sqrt{\ln k_n}$,
where the Lipschitz constant $L_n(C)$ is subject to some growth
restriction [see assumption (A4)].
\item[(A4)] 
%
%
\begin{equation} \label{eqhyp-26}
\forall C>0\qquad L_n(C) \delta_n k_n \sqrt{\ln k_n} \rightarrow0
\quad\mbox{and}\quad \eta_n k_n^{3/2}\sqrt{\ln k_n} \rightarrow0.
\end{equation}
Further, the design satisfies a condition on the trace of $\Phi^T\Phi$:
%
%
\begin{equation} \label{eqhyp-27}
\tr(\Phi^T\Phi) = O(n k_n).
\end{equation}
\end{longlist}

Since $\covmat$ is an orthogonal projection matrix on a
$k_n$-dimensional space, $\tr(\covmat)=k_n$ and $\eta_n^2 \geq k_n/n$.
Thus, the last part of (\ref{eqhyp-26}) implies $k_n^4 \ln k_n = o(n)$.

If we add the normality of the errors and a slight technical condition
$\ln n = o(k_n \ln k_n)$, these assumptions imply ours. Indeed, set
$M_n = C^2 k_n^2 \ln k_n$ for some arbitrary value of $C$. Our
condition~\ref{hypk-n-M-n} is immediate. Condition~\ref
{hypcontinuous-w} is got from (\ref{eqhyp-24}) and the first part of
(\ref{eqhyp-26}). The beginning\vspace*{1pt} of (A3) implies $-\ln w(\theta_0) =
O(k_n) = o(M_n)$. Using the concavity of the $\ln$ function and (\ref
{eqhyp-27}), we get $\ln\det(\Phi^T\Phi) \leq k_n \ln\tr(\Phi
^T\Phi
) - k_n \ln k_n = O(k_n \ln n) = o(M_n)$. Therefore, our condition \ref
{hypweight} holds.

\section{Semiparametric Bernstein--von Mises theorems} \label
{secsemiparametric}

We consider two kinds of functionals of $F$: linear and nonlinear ones.
These results can be easily adapted to functionals of $\theta$, using
the maps $\theta\mapsto\Phi\theta$ and $F\mapsto(\Phi^T\Phi
)^{-1} \Phi
^T F$.

\subsection{The linear case} \label{seclinear}

For linear functionals of $F$, we have the following corollary:
\begin{corol} \label{corollinear}
Let $p\geq1$ be fixed, and $G$ be a $\R^p \times\R^{n}$-matrix.
Suppose that the conditions of either Theorems \ref
{thmall-Gaussian-BvM} or \ref{thmsmooth-gaussian-BvM} are
satisfied. Then
\[
\E\bigl\| W(\ud(G F) | Y) - \mathcal{N}\bigl(G Y_{\phispan}, \sigma
_n^2 G \covmat G^T\bigr) \bigr\|_{\mathrm{TV}} \rightarrow0
\qquad\mbox{as } n\rightarrow\infty.
\]
Further, the distribution of $G Y_{\phispan}$ is $\mathcal{N}(G
F_{\phispan}, \sigma_n^2 G \covmat G^T)$.
\end{corol}

Corollary \ref{corollinear} is just a linear transform of the
preceding theorems, and of the distribution of $Y_{\phispan}$.

An example of application is given in Section \ref{secFourier}, in
the context of the regression on Fourier's basis.

\subsection{The nonlinear case} \label{secnonlinear}

The Bernstein--von Mises theorem which is presented here for nonlinear
functionals is derived from the nonparametric theorems thanks to Taylor
expansions. In the Taylor expansion of a functional, the first order
term naturally leads to the posterior normality, as in the case of
linear functionals. We do not want that the second order term interfere
with this phenomenon: it has to be controlled. The conditions of
Theorem~\ref{thmsmooth-gaussian-BvM} below are stated to permit this
control of the second order term.\looseness=1

Let $p\geq1$ be fixed, and $G\dvtx\R^{n} \mapsto\R^p$ be a twice
continuously differentiable function. For $F\in\R^{n}$, let $\dot{G}_F$
denote the Jacobian matrix of $G$ at $F$, and $D_F^2 G(\cdot, \cdot)$
the second derivative of $G$, as a bilinear function on $\R^{n}$. For
any $F\in\phispan$ and $a>0$, let
%
%
\begin{equation} \label{eqdef-B-a}
B_{F}(a) = \sup_{h\in\phispan\dvtx\|h\|^2 \leq\sigma_n^2 a} \sup
_{0\leq
t\leq1} \| D_{F+t h}^2 G(h, h)\|,
\end{equation}
where \mbox{$\|\cdot\|$} denotes the Euclidean norm of $\R^p$.

We also consider the following nonnegative symmetric matrix
%
%
\begin{equation} \label{eqGamma}
\Gamma_F = \sigma_n^2 \dot{G}_F \covmat\dot{G}_F^T.
\end{equation}
In the following, $\| \Gamma_F^{-1} \|$ denotes the Euclidean operator
norm of $\Gamma_F^{-1}$, which is also the inverse of the smallest
eigenvalue of $\Gamma_F$.


Let $\mathcal{I}$ be the collection of all intervals in $\R$, and for
any $I \in\mathcal{I}$, let $\psi(I)= \Proba(Z\in I)$, where $Z$ is a
$\mathcal{N}(0,1)$ random variable. Recall that $Y_\phispan$ is the MLE
and the orthogonal projection of $Y$ on $\phispan$.
\begin{theorem} \label{thmfunctional-base}
Let $G\dvtx\R^{n} \mapsto\R^p$ be a twice continuously differentiable
function, and let $\Gamma_F$ be as just defined. Suppose that $\Gamma
_{F_\phispan}$ is nonsingular, and that there exists a sequence
$(M_n)_{n\geq1}$ such that $k_n = o(M_n)$ and
%
%
\begin{equation} \label{eqhyp-B}
B_{F_\phispan}^2(M_n) = o\bigl(\bigl\|\Gamma_{F_\phispan}^{-1}\bigr\|
^{-1}\bigr).
\end{equation}
Suppose further that the conditions of either Theorems \ref
{thmall-Gaussian-BvM} or \ref{thmsmooth-gaussian-BvM} are
satisfied. Then, for any $b\in\R^p$,
%
%
\begin{equation} \label{eqbayesian-functionals}\qquad
\E\biggl[\sup_{I\in\mathcal{I}} \biggl| W\biggl( \frac{b^T
(G(F)-G(Y_\phispan))}{\sqrt{b^T\Gamma_{F_\phispan} b}} \in I
\Big| Y\biggr) - \psi(I) \biggr|\biggr] \rightarrow0
\qquad\mbox{as }
n\rightarrow\infty.
\end{equation}
Under the same conditions,
%
%
\begin{equation} \label{eqfrequentist-functionals}
\sup_{I\in\mathcal{I}} \biggl| \Proba\biggl( \frac{b^T (G(Y_\phispan
)-G(F_\phispan))}{\sqrt{b^T\Gamma_{F_\phispan} b}} \in I \biggr)
- \psi(I) \biggr| \rightarrow0 \qquad\mbox{as } n\rightarrow\infty.
\end{equation}
\end{theorem}

Note that $\sup_{I\in\mathcal{I}} |Q(I) - Q'(I)|$ is the
Levy--Prokhorov distance between two distributions $Q$ and $Q'$ on $\R$.
The Levy--Prokhorov distance metrizes the convergence in distribution.
So, when $p=1$ the Levy--Prokhorov distance between the distribution
$W(\ud G(F)| Y)$ and $\mathcal{N}(G(Y_\phispan),\Gamma
_{F_\phispan})$ goes to $0$ in mean, while $G(Y_\phispan
)$ goes to $\mathcal{N}(G(F_\phispan),\Gamma
_{F_\phispan})$ in distribution.

An application of Theorem \ref{thmfunctional-base} is given in
Section \ref{secFourier}, in the context of the regression on
Fourier's basis. The proof is delayed to Section \ref
{secproof-functional-all-gaussian}.

\section{Applications} \label{secapplications}

Here we give the three applications described in Section~\ref
{secFramework}. The models studied and the collections of regressors
used have already been defined there.

\subsection{The Gaussian sequence model} \label{secwhitenoise}

We consider the model (\ref{eqwhite-noise}). Here the MLE is the
projection $\theta_Y = (Y_j)_{1\leq j\leq k_n}$.

The nonparametric\vspace*{1pt} case corresponds to the estimation of
$\theta^0$. Under the assumption that $\theta^0$ is in some regularity
class, we will obtain a Bernstein--von Mises theorem with the posterior
convergence rate already obtained in previous works, in particular,
\citet{GhosalVanderVaart07noniid}. On the other hand, for some
priors known to achieve this rate, it will be seen that the centering
point and the asymptotic variance of the posterior distribution do not
fit
with the ones expected in a Bernstein--von Mises theorem. 
We also look at the semiparametric estimation of the squared $\ell^2$
norm of~$\theta^0$.\looseness=1

\subsubsection{The nonparametric estimation of
$\theta^0$} \label{secGaussian-seq-nonparam}

%
\begin{prop} \label{propwhite-noise-all-gaussian}
Suppose that $\sum_{j=1}^{k_n} (\theta_j^0)^2$ is bounded. This holds
when~$\theta^0$ is an element of $\ell^2(\N)$ not depending on $n$.
With a prior $\widetilde{W}=\mathcal{N}(0, \tau_n^2 I_{k_n}
)$ such that $n^{-1/4} = o(\tau_n)$, we have for any sequence $k_n\leq n$,
\[
\E\biggl\| \widetilde{W}(\ud\theta| Y) - \mathcal{N}\biggl(\theta_Y,
\frac{1}{n} I_{k_n} \biggr) \biggr\|_{\mathrm{TV}} \rightarrow0
\qquad\mbox{as } n\rightarrow\infty,
\]
and the convergence rate of $\theta$ \textit{toward $\theta_0$} is
$\sqrt
{\frac{k_n}{n}}$: for every $\lambda_n\rightarrow\infty$,
\[
\E\Biggl[ \widetilde{W}\Biggl( \|\theta-\theta_0\| \geq\lambda_n
\sqrt{\frac{k_n}{n}} \Big| Y\Biggr)\Biggr] \rightarrow0.
\]
\end{prop}

Recall that $\theta_0 = (\theta_j^0)_{1\leq j\leq k_n}$ is the
projection of $\theta^0$.
\begin{pf*}{Proof of Proposition \ref{propwhite-noise-all-gaussian}}
The beginning is an immediate corollary of Theorem \ref
{thmall-Gaussian-BvM}. For the convergence rate, let $\lambda
_n\rightarrow\infty$. Since $\theta_Y-\theta_0\sim\mathcal{N}(0,
\frac{1}{n} I_{k_n})$,
\[
\Proba\Biggl( \|\theta_Y-\theta_0\| \geq\frac{\lambda_n}{2} \sqrt
{\frac
{k_n}{n}} \Biggr) \rightarrow0.
\]
In the same way
\begin{eqnarray*}
\E\Biggl[ \widetilde{W}\Biggl( \|\theta-\theta_Y\| \geq\frac{\lambda
_n}{2} \sqrt{\frac{k_n}{n}} \Biggr) \Biggr] &\leq&\E\biggl\| \widetilde
{W}(\ud\theta| Y) - \mathcal{N}\biggl(\theta_Y, \frac{1}{n}
I_{k_n}\biggr)\biggr \|_{\mathrm{TV}} \\
&&{} + \mathcal{N}\biggl(0, \frac{1}{n} I_{k_n}\biggr)\Biggl(\Biggl\{
\|h\|\leq\frac{\lambda_n}{2} \sqrt{\frac{k_n}{n}}\Biggr\} \Biggr),
\end{eqnarray*}
which goes to $0$. Therefore,
\[
\E\Biggl[ \widetilde{W}\Biggl( \|\theta-\theta_0\| \geq\lambda_n \sqrt
{\frac{k_n}{n}} \Biggr) \Biggr] \rightarrow0.
\]
\upqed\end{pf*}

However, in such a general setting we have no information about the
bias between $\theta^0$ and its projection $\theta_0$. Several authors
add the assumption that the true parameter belongs to a Sobolev class
of regularity $\alpha>0$, defined by the relation ${\sum_{j=1}^{\infty}}
|\theta^0_j|^2 j^{2\alpha} <\infty$. In this setting we 
show that for some priors the induced posterior may achieve the
nonparametric convergence rate but with a centering point and a
variance different from what is expected in the Bernstein--von Mises
theorem. Then we exhibit priors for which both the Bernstein--von Mises
theorem and the nonparametric convergence rate hold.

From now on,\vspace*{-1pt} we suppose that $\sum_{j=1}^{\infty} |\theta^0_j|^2
j^{2\alpha} <\infty$. In this setting \citet{GhosalVanderVaart07noniid},
Section 7.6, consider a prior $\widetilde{W}$ such that $\theta_1$,
$\theta_2,\ldots$ are independent, and $\theta_j$ is normally
distributed with variance $\sigma_{j,k_n}^2$. Further, the variances
are supposed to satisfy
%
%
\begin{equation} \label{eqprior-white-noise}
c/k_n \leq\min\{ \sigma_{j,k_n}^2 j^{2\alpha}\dvtx1\leq j\leq k_n \}
\leq C/k_n
\end{equation}
for some positive constants $c$ and $C$. Suppose that $\alpha\geq1/2$
and there exist constants $C_1$ and $C_2$ such that $C_1
n^{1/(1+2\alpha
)} \leq k_n \leq C_2 n^{1/(1+2\alpha)}$. Then \citet
{GhosalVanderVaart07noniid}, Theorem 11, proved that the posterior
converges at the rate $n^{-\alpha/(1+2\alpha)}$.

In order to get $n^{-1} I_{k_n}$ as asymptotic variance, we need more
stringent conditions on $k_n$, or a flatter prior. To see this is
necessary, consider, for $k_n \approx n^{1/(1+2\alpha)}$, the following
choice for $\sigma_{j, k_n}$:
\[
\sigma_{j, k_n}^2 = \cases{
k_n^{-1}, &\quad if $1 \leq j \leq k_n/2$,\vspace*{2pt}\cr
2^{2\alpha}/n, &\quad if $j > k_n/2$.}
\]
Then $\min\{ \sigma_{j,k_n}^2 j^{2\alpha}\dvtx1\leq j\leq k_n \}
\approx
k_n^{-1}$, and the posterior converges at the rate $n^{-\alpha
/(1+2\alpha)}$.

For this case we can explicitly calculate the posterior distribution.
This is similar to the calculation made in the proof of Theorem \ref
{thmall-Gaussian-BvM}.
The coordinates are independent, and
\[
\widetilde{W}(\ud\theta_j | Y) = \mathcal{N} \biggl(\frac{\sigma_{j,
k_n}^2}{\sigma_n^2 + \sigma_{j, k_n}^2} Y_j, \frac{\sigma_n^2
\sigma
_{j, k_n}^2}{\sigma_n^2 + \sigma_{j, k_n}^2} \biggr).
\]
For $j> k_n/2$, $\frac{\sigma_{j, k_n}^2}{\sigma_n^2 + \sigma_{j,
k_n}^2} = \frac{4^\alpha}{1+4^\alpha}$, and, therefore, $\|
\widetilde{W}(\ud\theta_j | Y) - \mathcal{N} (Y_j, \sigma_n^2
) \|_{\mathrm{TV}}$ is bounded away from $0$.

By contrast, with an isotropic and flat prior we obtain the centering
point and the asymptotic variance we expected, and the same convergence
rate as previously. We have the following:
\begin{prop} \label{propwhite-noise-all-gaussian-bias}
Suppose\vspace*{2pt} that $\theta^0$ belongs to the Sobolev class of
regularity $\alpha>0$. Choose a prior $\widetilde{W}=\mathcal{N}(0,
\tau_n^2 I_{k_n})$ such that $n^{-1/4} = o(\tau_n)$, which ensures the
asymptotic normality of the posterior distribution as in Proposition
\ref{propwhite-noise-all-gaussian}.\vspace*{1pt}

If further\vspace*{1pt} $k_n\approx n^{1/(1+2\alpha)}$, then the convergence rate
of $\theta$ toward $\theta_0$ and toward $\theta^0$ is $n^{-\alpha
/(1+2\alpha)}$: for every $\lambda_n\rightarrow\infty$,
\[
\E\bigl[ \widetilde{W}\bigl( \|\theta-\theta^0\| \geq\lambda_n
n^{-\alpha/(1+2\alpha)} | Y\bigr)\bigr] \rightarrow0.
\]
\end{prop}
\begin{pf}
We consider $\theta$ and $\theta_0$ as elements of $\ell^2(\N)$ by
setting $\theta_j = \theta_{0,j}=0$ for $j\geq k_n+1$. The convergence
rate toward $\theta_0$ has already been established in Proposition
\ref
{propwhite-noise-all-gaussian}. Since $\theta_{0,j} = \theta^0_j$ for
$1\leq j\leq k_n$, $\|\theta^0-\theta_0\| \leq k_n^{-\alpha} \sqrt
{\sum
_{j=k_n+1}^{\infty} (\theta^0_j)^2 j^{2\alpha} } = O(k_n^{-\alpha
})$. Therefore, the convergence rate of $\theta$ toward $\theta
^0$ is also $n^{-\alpha/(1+2\alpha)}$.
\end{pf}

\subsubsection{Semiparametric theorem for the
$\ell^2$ norm of $\theta^0$}
\label{secGaussian-seq-semiparam}

We consider the prior distribution 
used in Proposition \ref{propwhite-noise-all-gaussian-bias},
but\vspace*{1pt} now we look at the posterior distribution
of~$\|\theta\|^2$. To get asymptotic normality with variance
$n^{-1/2}$, we just need $k_n=o(\sqrt {n})$. To control the bias term,
we need $\alpha>1/2$, and in this case we get an adaptive Bayesian
estimator.
\begin{prop} \label{propwhite-noise-l2-norm}
Let $\alpha>1/2$ and suppose that $\theta^0$ belongs to the Sobolev
class of regularity $\alpha$. Choose a prior $\widetilde{W}=\mathcal
{N}(0, \tau_n^2 I_{k_n})$ such that $n^{-1/4} = o(\tau_n)$.
Then, for any choice of $k_n$ such that $k_n = o(\sqrt{n})$ and $\sqrt
{n} = o(k_n^{2\alpha})$,
\[
\E\biggl[\sup_{I\in\mathcal{I}} \biggl| \widetilde{W}\biggl( \frac
{\sqrt{n} (\|\theta\|^2 - \|\theta_Y\|^2)}{2 \|\theta^0\|}
\in I \Big| Y\biggr) - \psi(I) \biggr|\biggr] \rightarrow0
\qquad\mbox{as
} n\rightarrow\infty
\]
and\vspace*{1pt}
$\frac{\sqrt{n} (\|\theta_Y\|^2 - \|\theta_0\|^2 )}{2 \|
\theta^0\|} \rightarrow\mathcal{N}(0,1)
$ in distribution, as $n\rightarrow\infty$. Further, the bias is
negligible with respect to the square root of the variance:
\[
\frac{\sqrt{n} (\|\theta_0\|^2 - \|\theta^0\|^2 )}{2 \|\theta
^0\|} = o(1).
\]
In particular, the choice $k_n = \sqrt{n/\ln n}$ is adaptive in
$\alpha$.
\end{prop}
\begin{pf}
We set up an application of Theorem \ref{thmfunctional-base}. Since
$\sigma_n=n^{-1/2}$, the conditions of Theorem \ref
{thmall-Gaussian-BvM} are fulfilled. 

Here $G(\theta) = \theta^T \theta$, $\dot{G}_\theta= 2 \theta^T$ and
$\ddot{G}_\theta= 2 I_{k_n}$. Therefore, $B_{\theta_0}(M_n) = 2
M_n/n$, while $\Gamma_{\theta_0} = 4 \|\theta_0\|^2/n$.

Let us choose $(M_n)_{n\geq1}$ such that $k_n=o(M_n)$ and $M_n=o(\sqrt
{n})$. Such sequences exist and fulfill the conditions of Theorem \ref
{thmfunctional-base}.

Since $\|\theta_0\|^2 \rightarrow\|\theta^0\|^2$, we can substitute
the variance $\Gamma_{\theta_0}$ by $4 \|\theta^0\|^2/n$ and get the
two asymptotic normality results, (\ref{eqbayesian-functionals}) and
(\ref{eqfrequentist-functionals}).

As\vspace*{1pt} $n\rightarrow\infty$, $\|\theta^0\|^2 - \|\theta_0\|^2 = \|
\theta
^0-\theta_0\|^2 = O(k_n^{-2\alpha})$, as in the proof of
Proposition \ref{propwhite-noise-all-gaussian-bias}. If $\sqrt{n} =
o(k_n^{2\alpha})$, we get $\sqrt{n} (\|\theta_0\|^2 - \|\theta^0\|
^2 ) = o(1)$.
\end{pf}

\subsection{Regression on Fourier's basis} \label{secFourier}

Now we consider the regression mo\-del~(\ref{eqregression-0-1}) with a
function $f$ in a Sobolev class $\mathcal{W}(\alpha, L)$, and use
Fourier's basis~(\ref{eqFourier-basis}). For any $\theta\in\R^{k_n}$,
we define $f_\theta= \sum_{j=1}^{k_n} \theta_j \varphi_j$. We also
denote by $\theta^0\in\ell^2(\N)$ the sequence of Fourier's
coefficients of $f$: $f=\sum_{j=1}^\infty\theta^0_j \varphi_j$.

The following useful lemma about our collection of regressors can be
found, for instance, in \citet{Tsybakov04} (we slightly modified
it to take into account the case~$n$ even):
\begin{lem} \label{lemTsibakov}
Suppose either that $n$ is odd and $k_n\leq n$, or $n$ is even and
$k_n\leq n-1$. Consider the collection $(\phi_j)_{1\leq j\leq k_n}$
defined before, and $\Phi$ the associated matrix. Then
\[
\Phi^T \Phi= n I_{k_n}.
\]
\end{lem}

This makes the regression on Fourier's basis very close to the Gaussian
sequence model, and the results we obtain are similar.

In this subsection we first consider the estimation of $f$ in a Sobolev
class, for which we get a Bernstein--von Mises theorem and the
frequentist minimax $n^{-\alpha/(1+2\alpha)}$ posterior convergence
rate for the $L^2$ norm. Then we consider two semiparametric settings:
the estimation of a linear functional of $f$, and the estimation of the
$L^2$ norm of $f$. We get the adaptive $\sqrt{n}$ convergence rate for
any $\alpha>1/2$.

\subsubsection{Nonparametric Bernstein--von Mises theorem in Sobolev
classes} \label{secSobolev-nonparam}

\begin{prop} \label{propfourier-nonparametric}
Suppose that $f$ belongs to some Sobolev class $\mathcal{W}(\alpha,
L)$ for $L>0$ and $\alpha>1/2$. Let $k_n\approx n^{1/(1+2\alpha)}$ and
$\widetilde{W}=\mathcal{N}(0, \gamma_n I_{k_n})$ be the
prior on 
$\theta$, for a sequence $(\gamma_n)_{n\geq1}$ such that $1/\sqrt{n}
= o(\gamma_n)$. Then
\[
\E\biggl\| \widetilde{W}(\ud\theta| Y) - \mathcal{N}\biggl(\theta_Y,
\frac{\sigma^2}{n} I_{k_n}\biggr) \biggr\|_{\mathrm{TV}} \rightarrow0
\qquad\mbox{as } n\rightarrow\infty
\]
and the convergence rate relative to the Euclidean norm for $f_\theta$
is $n^{-\alpha/(1+2\alpha)}$: for every $\lambda_n\rightarrow\infty$,
\[
\E\bigl[ \widetilde{W}\bigl( \|f_\theta-f\| \geq\lambda_n
n^{-\alpha/(1+2\alpha)} | Y\bigr)\bigr] \rightarrow0.
\]
\end{prop}
\begin{pf}
The conditions of Theorem \ref{thmall-Gaussian-BvM} are fulfilled:
with $\tau_n^2 = n \gamma_n$, we have $n=o(\tau_n^4)$. The first
assertion follows.

Because of the orthogonal nature of Fourier's basis, $\|f_\theta-f\| =
\|\theta-\theta^0\|$ in $\ell^2(\N)$. We use the decomposition $\|
\theta
-\theta^0\|^2 \leq\|\theta-\theta_0\|^2 + \|\theta_0-\theta^0\|
^2$. In
the same way as in the proof of Proposition \ref
{propwhite-noise-all-gaussian}, for any $\lambda_n\rightarrow\infty$,
\[
\E\Biggl[ \widetilde{W}\Biggl( \|\theta-\theta_0\| \geq\lambda_n \sqrt
{\frac{k_n}{n}} \Biggr) \Biggr] \rightarrow0.
\]
Going back to Definition \ref{dfSobolev}, we have
\[
\|\theta_0-\theta^0\|^2 = \sum_{j=k_n+1}^\infty(\theta_j^0)^2 \leq
k_n^{-2\alpha} \sum_{j=k_n+1}^\infty a_j^{2\alpha} (\theta_j^0)^2 =
O(k_n^{-2\alpha}).
\]
This permits to get
\[
\E\bigl[ \widetilde{W}\bigl( \|\theta-\theta^0\| \geq\lambda_n
n^{-\alpha/(1+2\alpha)} | Y\bigr)\bigr] \rightarrow0.
\]
\upqed\end{pf}

\subsubsection{Linear functionals of $f$}
\label{seclin-func-of-f}

Let $g\dvtx[0, 1] \rightarrow\R$ be\vspace*{1pt} a function in $\LL^2([0,1])$. We want
to estimate $\mathcal{F}(f)= \int_0^1 f g$, and we approximate it by
\[
\frac{1}{n} \sum_{i=1}^n g(i/n) f(i/n) = G F_0,
\]
where\vspace*{1pt} $G = (g(i/n)/n)_{1\leq i\leq n}^T$. The plug-in MLE
estimator of $G F_0$ in the misspecified model $\phispan$ is $G
Y_\phispan$. More generally, we consider the functional $F\mapsto G F$.
The following result is adaptive, in the sense that the same choice
$k_n =\lfloor n/\ln n \rfloor$ entails the convergence rate
$n^{-1/2}$ for all values of $\alpha>1/2$.
\begin{prop} \label{corollinear-functional-of-f}
Suppose $f$ is bounded, and let $W$ be the prior induced by the
$\mathcal{N}(0, \gamma_n$ $I_{k_n})$ distribution on $\theta$, for a
sequence $(\gamma_n)_{n\geq1}$ such that $1/\sqrt{n} = o(\gamma
_n)$. Then:
\renewcommand\thelonglist{(\arabic{longlist})}
\renewcommand\labellonglist{\thelonglist}
\begin{longlist}
\item\label{lin-func-no-hyp}
\[
\E\bigl\| W(\ud(G F) | Y) - \mathcal{N}\bigl(G Y_{\phispan}, \sigma
^2 G \covmat G^T\bigr) \bigr\|_{\mathrm{TV}} \rightarrow0
\]
and the distribution of $G Y_{\phispan}$ is $\mathcal{N}(G
F_{\phispan}, \sigma^2 G \covmat G^T)$.
\item\label{lin-func-sobolev}
Suppose further\vspace*{1pt} that $f$ and $g$ belong to some Sobolev class $\mathcal
{W}(\alpha, L)$ for $L>0$ and $\alpha>1/2$. Then $G \covmat G^T \sim
\frac{1}{n} \int_0^1 g^2$,
\[
\E\biggl\| W \biggl( \ud\frac{\sqrt{n}(G F - G Y_{\phispan
})}{\sigma\sqrt{\int_0^1 g^2}} \Big| Y\biggr) - \mathcal{N}(0,
1) \biggr\|_{\mathrm{TV}} \rightarrow0
\]
and
$
\frac{\sqrt{n}(G Y_{\phispan} - G F_\phispan)}{\sigma\sqrt{\int_0^1
g^2}} \rightarrow\mathcal{N}(0,1)
$ in distribution, as $n \rightarrow\infty$.
\item\label{lin-func-bias}
Suppose that $f$ and $g$ belong to some Sobolev class $\mathcal
{W}(\alpha, L)$ for $L>0$ and $\alpha>1/2$, and suppose further that
$k_n$ is large enough so that $n=o(k_n^{2\alpha})$. Then the bias is
negligible with respect to the square root of the variance:
\[
\frac{\sqrt{n} (G F_\phispan- \mathcal{F}(f))}{\sigma\sqrt
{\int_0^1 g^2}} = o(1).
\]
\end{longlist}
\end{prop}

Before the proof we give two lemmas, proved in Appendix B %
in the supplemental article [\citet{BvMGRegSupplement}], about the error
terms of the approximation of a Sobolev class by a sieve build on
Fourier's basis, and of the approximation of an integral by a Riemann sum.
\begin{lem} \label{lemprojection-bias}
Let $\alpha>1/2$ and $L>0$. We suppose $n$ odd or $k_n<n$. If $f\in
\mathcal{W}(\alpha, L)$,
\[
\bigl\| F_0- F_\phispan\bigr\| \leq\bigl(1+o(1)\bigr) \frac{\sqrt{2}L}{\pi^\alpha}
\frac
{\sqrt{n}}{k_n^\alpha}.
\]
Further, $\|F_0\| \sim\sqrt{n \int_0^1 f^2}$ and $\|F_0-F_\phispan
\| =
O(k_n^{-\alpha} \|F_0\|)$.
\end{lem}
\begin{lem} \label{lemRiemann-bias}
Let two functions $f\in\mathcal{W}(\alpha, L)$ and $g\in\mathcal
{W}(\alpha', L')$ for some $\alpha, \alpha'>1/2$ and two positive
numbers $L$ and $L'$. Then
\[
\Biggl| \frac{1}{n} \sum_{i=1}^n f(i/n) g(i/n) -\int_0^1 f g \Biggr| =
O\bigl(n^{-\inf(\alpha, \alpha')}\bigr).
\]
\end{lem}
\begin{pf*}{Proof of Proposition \ref{corollinear-functional-of-f}}
(1) The first assertion is just Corollary \ref{corollinear}. The
conditions of Theorem \ref{thmall-Gaussian-BvM} are fulfilled, as in
the proof of Proposition \ref{propfourier-nonparametric}.\vadjust{\goodbreak}

(2) If\vspace*{1pt} $g \in\mathcal{W}(\alpha, L)$ for $L>0$ and $\alpha>1/2$, $G
\covmat G^T = \| \covmat G^T\|^2 \sim\| G^T\|^2$ by Lemma \ref
{lemprojection-bias}. In the meantime $\| G^T\|^2 = \frac{1}{n^2}
\sum
_{i=1}^n g^2(x_i) \sim\frac{1}{n} \int_0^1 g^2$ by Lemma \ref
{lemRiemann-bias}. So $G \covmat G^T \sim\frac{1}{n} \int_0^1 g^2$,
and the variance in the formulas of Corollary \ref{corollinear} can be
substituted with $\frac{1}{n} \int_0^1 g^2$.

(3) We decompose the bias into two terms, $|G F_0 - \mathcal{F}(f)|$
and $|G F_\phispan- G F_0|$, and show that both are $o(n^{-1/2})$. The
first term is controlled by Lemma \ref{lemRiemann-bias}. For the last
one, $|G F_\phispan- G F_0| \leq\|G^T\| \|F_\phispan- F_0\|$.
But\vspace*{1pt} $\|G^T\|=O(n^{-1/2})$, $\|F_\phispan- F_0\| =
O(k_n^{-\alpha} \|F_0\| )$ by Lemma \ref{lemprojection-bias} and
$\|F_0\| = O(\sqrt{n})$. We conclude thanks to the assumption
$n=o(k_n^{2\alpha})$.
\end{pf*}

\subsubsection{$L^2$ norm of $f$}
\label{secL2}

Suppose that we want to estimate $\mathcal{F}(f) = \int_0^1 f^2$. We
can consider the plug-in MLE estimator
\[
G\bigl(Y_\phispan\bigr) = \frac{1}{n} \bigl\|Y_\phispan\bigr\|^2 =
\frac{1}{n} \sum_{i=1}^{n} \Biggl( \sum_{j=1}^{k_n} \theta_{Y, j}
\varphi_j(i/n) \Biggr)^2.
\]
More generally, we define, for any $F \in\R^{n}$,
%
%
\begin{equation} \label{eqG-f2}
G(F) = \frac{1}{n} \| F \|^2. 
\end{equation}

With a Gaussian prior, we obtain the following result, which is also
adaptive: the same $k_n=\lfloor\sqrt{n}/\ln n\rfloor$ is suitable
whatever $\alpha>1/2$.
\begin{prop} \label{corolL2-norm}
Let $G(F)=\|F\|^2/n$. Suppose that $f\in\mathcal{W}(\alpha, L)$ for
some $L>0$ and $\alpha>1/2$. Let $W$ be the prior induced by the
$\mathcal{N}(0, \gamma_n$ $I_{k_n})$ distribution on $\theta$, for a
sequence $(\gamma_n)_{n\geq1}$ such that $1/\sqrt{n} = o(\gamma_n)$.
The sequence $(k_n)_{n\geq1}$ can be chosen such that $k_n = o(\sqrt
{n})$ and $\sqrt{n} = o(k_n^{2\alpha})$, and with such a choice,
\[
\E\biggl[\sup_{I\in\mathcal{I}} \biggl| W\biggl( \frac{\sqrt{n}
(G(F)-G(Y_\phispan))}{2\sigma\sqrt{\mathcal{F}(f)}} \in I
\Big| Y\biggr) - \psi(I) \biggr|\biggr] \rightarrow0
\qquad\mbox{as }
n\rightarrow\infty
\]
and
$
\frac{\sqrt{n} (G(Y_\phispan)-G(F_\phispan))}{2\sigma\sqrt
{\mathcal{F}(f)}} \rightarrow\mathcal{N}(0,1)
$ in distribution, as $n\rightarrow\infty$. Further, the bias is
negligible with respect to the square root of the variance:
\[
\frac{\sqrt{n} (G(F_\phispan) - \mathcal{F}(f))}{2\sigma
\sqrt{\mathcal{F}(f)}} = o(1).
\]
\end{prop}

A similar corollary could be stated for a non-Gaussian prior.
\begin{pf*}{Proof of Proposition \ref{corolL2-norm}}
First, let us note that the conditions of Theorem~\ref
{thmall-Gaussian-BvM} are fulfilled, as in the proof of
Proposition \ref{propfourier-nonparametric}. Lemma 10 %
in Appendix~B 
insures that $f$ is bounded.

In this setting $\dot{G}_F = (2/n) F^T$ and $D_F^2 G(h, h) = (2/n)
\|h\|^2$ for any $F\in\R^n$ and any $h\in\R^{n}$. Therefore, $B_{F}(a)
= 2 \sigma^2 a/n$, and $\Gamma_F = 4 (\sigma^2 / n^2) \|F\|^2$. By
Lem\-ma~\ref{lemprojection-bias}, $\|F_\phispan\|^2 \sim\|F_0\|^2
\sim
n \mathcal{F}(f)$. Thus, $\Gamma_{F_\phispan} = 4 (1+o(1)) \mathcal
{F}(f)/n$.

Let us choose $(M_n)_{n\geq1}$ such that $k_n=o(M_n)$ and $M_n=o(\sqrt
{n})$. Such sequences exist and fulfill the conditions of Theorem \ref
{thmfunctional-base}. We can substitute the variance $\Gamma
_{F_\phispan}$ by $4 \mathcal{F}(f)/n$ and get the two asymptotic
normality results.

Let us now consider the bias term:
\[
\mathcal{F}(f) - G\bigl(F_\phispan\bigr) \leq\frac{\|F_0\|^2 -
\|F_\phispan\| ^2}{n} + \Biggl(\int_0^1 f^2 - \frac{1}{n}\sum_{i=1}^n
f^2(i/n) \Biggr).
\]
We use Lemma \ref{lemprojection-bias} to control $\|F_0\|^2 - \|
F_\phispan\|^2$, and Lemma \ref{lemRiemann-bias} for the other term:
\[
\bigl| \mathcal{F}(f) - G\bigl(F_\phispan\bigr) \bigr| = O(k_n^{-2\alpha
}) + O(n^{-\alpha}).
\]
This is a $o(1/\sqrt{n})$ under the assumptions of Corollary \ref
{corolL2-norm}.
\end{pf*}

\subsection{Regression on splines} \label{secsplines}

Here we consider the regression model for functions in $C^\alpha[0,1]$
with $\alpha>0$, using splines, set up in Section \ref{secFramework}.
We first develop further the framework and the assumptions used here,
and recall the previous result of \citet{GhosalVanderVaart07noniid},
Section 7.7.1, which obtains the posterior concentration at the
frequentist minimax rate. Then we present two Bernstein--von Mises
theorems: the first one with the same prior as \citet
{GhosalVanderVaart07noniid} but a stronger condition on $k_n$ (or
equivalently on $\alpha$); the second one with a flatter prior, for
which we obtain the minimax convergence rate in addition to the
asymptotic Gaussianity of the posterior distribution.

To see this, we begin with some preliminaries. For any $\theta\in\R
^{k_n}$, define $f_\theta= \sum_{j=1}^{k_n} \theta_j B_j$. The
$B$-splines basis has the following approximation property: for any
$\alpha>0$, there exist $C_{\alpha}>0$ such that, if $f\in C^\alpha
[0,1]$, there exists $\theta^\infty\in\R^{k_n}$ satisfying
%
%
\begin{equation} \label{eqspline-approx}
\| f - f_{\theta^\infty} \|_{\infty} \leq C_{\alpha}
k_n^{-\alpha} \|f\|_{\alpha}.
\end{equation}

We need the design $(x_i^{(n)})_{n\geq1, 1\leq i \leq n}$ to
be sufficiently regular and, as stressed in \citet
{GhosalVanderVaart07noniid}, the spatial separation property of
$B$-splines permits us to express the precise condition in terms of the
covariance matrix $\Phi^T\Phi$. We suppose that there exist positive
constants $C_1$ and $C_2$ such that, as $n$ increases, for any $\theta
\in\R^{k_n}$,
%
%
\begin{equation} \label{eqspline-design-regularity}
C_1 \frac{n}{k_n} \|\theta\|^2 \leq\theta^T \Phi^T\Phi\theta\leq
C_2 \frac{n}{k_n} \|\theta\|^2.
\end{equation}

Let us associate the norm $\|f\|_n=\sqrt{{\frac{1}{n} \sum_{i=1}^{n}}
|f(x_i)|^2}$ to the design. Note that $\sqrt{n} \|f_{\theta}\|_n= \|
\Phi
\theta\|$ if $\theta\in\R^{k_n}$. Under (\ref
{eqspline-design-regularity}) we have a relation between \mbox{$\|\cdot\|_n$}
and the Euclidean norm on the parameter space: for every $\theta_1$
and~$\theta_2$,
\[
C_1 \|\theta_1-\theta_2\| \leq\sqrt{k_n} \|f_{\theta_1} -
f_{\theta_2}\|_n \leq C_2 \|\theta_1-\theta_2\|.
\]

With these conditions \citet{GhosalVanderVaart07noniid}, Theorem
12, get
the posterior concentration at the minimax rate. Take $\alpha\geq1/2$,
let $\widetilde{W}=\mathcal{N}(0, I_{k_n})$ be the prior on
the spline coefficients, and suppose there exist constants $C_3$ and
$C_4$ such that $C_3 n^{1/(1+2\alpha)} \leq k_n \leq C_4
n^{1/(1+2\alpha
)}$. Then the posterior concentrates at the minimax rate $n^{-\alpha
/(1+2\alpha)}$ relative to $\|\cdot\|_n$: for every $\lambda
_n\rightarrow\infty$,
\[
\E\bigl[ \widetilde{W}\bigl( \|f_\theta-f\|_n \geq\lambda_n
n^{-\alpha/(1+2\alpha)} | Y\bigr)\bigr] \rightarrow0.
\]
This is equivalent to a convergence rate $n^{({1-2\alpha
})/({2(1+2\alpha)})}$ relative to the Euclidean norm for~$\theta$:
\[
\E\bigl[ \widetilde{W}\bigl( \|\theta-\theta_0\|\geq\lambda_n
n^{({1-2\alpha})/({2(1+2\alpha)})} | Y\bigr)\bigr] \rightarrow0.
\]
Indeed, (\ref{eqspline-approx}) and the projection property imply
\[
\|f_{\theta_0} - f\|_n \leq\|f_{\theta^\infty} - f\|_n \leq\|
f_{\theta
^\infty} - f\|_\infty\leq C_{\alpha} \|f\|_{\alpha} k_n^{-\alpha}.
\]

Now, with modified assumptions we get the Bernstein--von Mises theorem
in two different settings. First, with the same prior as \citet
{GhosalVanderVaart07noniid}:\vspace*{-3pt}
\begin{prop} \label{propBvMsplineswiththeirprior}
Assume that $f$ is bounded, $k_n=o((\frac{n}{\ln n}
)^{1/3})$ and (\ref{eqspline-design-regularity}) holds. Let
$\widetilde{W}=\mathcal{N}(0, I_{k_n})$ be the prior on the
spline coefficients. Then
%
%
\begin{equation} \label{eqBvM-splines-1}
\E\bigl\| \widetilde{W}(\ud\theta| Y) - \mathcal{N}(\theta_Y,
\sigma^2 (\Phi^T \Phi)^{-1}) \bigr\|_{\mathrm{TV}} \rightarrow0
\qquad\mbox{as } n\rightarrow\infty
\end{equation}
and the convergence rate relative to the Euclidean norm for $\theta$ is
$\frac{k_n}{\sqrt{n}}$.\vspace*{-3pt}
\end{prop}

\subsubsection*{Remarks} We need $\alpha>1$ to get the Gaussian shape with
the same convergence rate as in \citet{GhosalVanderVaart07noniid}. The
conditions of Proposition \ref{propBvMsplineswiththeirprior} are
satisfied, in particular, if there\vspace*{1pt} exist constants $C_3$ and $C_4$ such
that $C_3 n^{1/(1+2\alpha)} \leq k_n \leq C_4 n^{1/(1+2\alpha)}$. In
this case the convergence rate for $\theta$ is $n^{({1-2\alpha
})/({2(1+2\alpha)})}$.
\begin{pf*}{Proof of Proposition \ref{propBvMsplineswiththeirprior}}
We set up an application of Theorem \ref{thmsmooth-gaussian-BvM}. We
can choose $M_n$ such that $k_n \ln n =o(M_n)$ and $M_n = o(\frac
{n}{k_n^2})$. Assumption \ref{hypk-n-M-n} is then trivially satisfied.

From (\ref{eqspline-design-regularity}) we get $\|\Phi^T\Phi
\| \leq C_2 \frac{n}{k_n}$ and $\|(\Phi^T\Phi)^{-1}\|
\leq C_1^{-1} \frac{k_n}{n}$. We have also $\ln\det(\Phi^T\Phi)
\leq
k_n \ln C_2 + k_n \ln(\frac{n}{k_n}) = O(k_n \ln n) =
o(M_n)$. Since $\theta_0 = \Phi(\Phi^T\Phi)^{-1} F_0$,
\[
\|\theta_0\|^2 \leq\frac{k_n}{C_1 n} \|F_0\|^2 \leq\frac{\|f\|
_\infty
}{C_1} k_n.
\]
Therefore, $-\ln w(\theta_0) = O(1) + \frac{1}{2} \|\theta_0\|^2 =
O(k_n) = o(M_n)$, and assumption~\ref{hypweight} holds.\vadjust{\goodbreak}

Let $h\,{\in}\,\R^{k_n}$ such that $\|\Phi h\|^2\,{\leq}\,\sigma^2 M_n$. We have
$\|h\|^2\,{\leq}\,\| (\Phi^T\Phi)^{-1} \|\cdot\break\|\Phi h\|^2
\,{\leq}\,\frac{\sigma^2 k_n M_n}{C_1 n}\,{=}\,o(k_n^{-1})$. Therefore,
%
%
\begin{equation} \label{eqto-get-cond-1}
\sup_{\|\Phi h\|^2\leq\sigma^2 M_n} \biggl|{\ln\frac{w(\theta
_0+h)}{w(\theta_0)} }\biggr| \leq\sup_{\|\Phi h\|^2\leq\sigma^2 M_n}
\frac{\|h\|^2 + 2 \|h\| \|\theta_0\|}{2} = o(1)
\end{equation}
and assumption \ref{hypcontinuous-w} follows.

Let us now prove the convergence rate. Let $\lambda_n\rightarrow
\infty
$. Then
\[
\Proba\biggl( \|\theta_Y-\theta_0\| \geq\frac{\lambda_n k_n}{2\sqrt
{n}} \biggr) \leq\Proba\biggl( \|\Phi(\theta_Y-\theta_0)\|^2 \geq\frac
{C_1 \lambda_n^2 k_n}{4} \biggr) \rightarrow0
\]
since $\|\Phi(\theta_Y-\theta_0)\|^2 \sim\sigma^2 \chi^2(k_n)$. In
the same way
\begin{eqnarray*}
\E\biggl[ \widetilde{W}\biggl( \|\theta-\theta_Y\| \geq\frac{\lambda_n
k_n}{2\sqrt{n}} \biggr) \biggr] &\leq&
\E\bigl\| \widetilde{W}(\ud
\theta| Y) - \mathcal{N}(\theta_Y, \sigma^2 (\Phi^T \Phi
)^{-1}) \bigr\|_{\mathrm{TV}} \\
&&{} + \mathcal{N}(0, \sigma^2 (\Phi^T \Phi)^{-1})
\biggl(\biggl\{h\dvtx\|h\| \leq\frac{\lambda_n k_n}{2\sqrt{n}}
\biggr\} \biggr)
\\
&\rightarrow&0,
\end{eqnarray*}
where Theorem \ref{thmsmooth-gaussian-BvM} controls the first term in
the right. Therefore, assumption~\ref{hypweight} holds:
\[
\E\biggl[ \widetilde{W}\biggl( \|\theta-\theta_0\| \geq\frac{\lambda_n
k_n}{\sqrt{n}} \biggr) \biggr] \rightarrow0.
\]
Now, (\ref{eqBvM-splines-1}) is the same as Theorem \ref
{thmsmooth-gaussian-BvM} in terms of $\widetilde{W}$.
\end{pf*}

The situation is similar to the one we encountered with the Gaussian
sequence model. To get the Bernstein--von Mises theorem with the same
convergence rate as \citet{GhosalVanderVaart07noniid} for $\alpha
\leq
1$, we need a~flatter prior: 
\begin{prop}\label{pr8}
Assume that $f$ is bounded and (\ref{eqspline-design-regularity})
holds. Let $\widetilde{W}=\mathcal{N}(0, \tau_n^2 I_{k_n})$
be the prior on the spline coefficients, with the sequence $\tau_n$ satisfying
\[
\frac{k_n^2 \ln n}{n} = o(\tau_n^2) \quad\mbox{and}\quad \frac
{k_n^3 \ln n}{n} = o(\tau_n^4).
\]
Then
\[
\E\bigl\| \widetilde{W}(\ud\theta| Y) - \mathcal{N}(\theta_Y,
\sigma^2 (\Phi^T \Phi)^{-1}) \bigr\|_{\mathrm{TV}} \rightarrow0
\qquad\mbox{as } n\rightarrow\infty
\]
and the convergence rate relative to the Euclidean norm for $\theta$ is
$\frac{k_n}{\sqrt{n}}$.
\end{prop}

When $\alpha>0$ and $k_n$ is of order $n^{1/(1+2\alpha)}$, the
conditions reduce to\break $n^{({2-2\alpha})/({1+2\alpha})} \ln n = o(\tau
_n^4)$. So we obtain the convergence rate\vadjust{\goodbreak} of \citet
{GhosalVanderVaart07noniid} in addition to the Gaussian shape with the
same $k_n$, even for $\alpha\leq1$, but with a different prior.\vspace*{-3pt}
\begin{pf*}{Proof of Proposition \ref{pr8}}
The proof is essentially the same as for Proposition \ref
{propBvMsplineswiththeirprior}. $M_n$ can be chosen so that $k_n
\ln n=o(M_n)$, $M_n=o(\frac{n\tau_n^2}{k_n})$, and
$M_n=o(\frac{n\tau_n^4}{k_n^2})$. These last two conditions
are the ones needed to obtain the same upper bounds as in (\ref
{eqto-get-cond-1}).
\end{pf*}

\vspace*{-3pt}\section{Proofs}\vspace*{-3pt} \label{secproofs}

\subsection{\texorpdfstring{Proof of Theorem \protect\ref{thmall-Gaussian-BvM}}{Proof of Theorem 1}}
\label{secproof-all-Gaussian}

In the present setting all distributions are explicit and admit known
densities with respect to the corresponding Lebesgue measure. We
decompose any $y\in\R^n$ in two orthogonal components $y = \Phi
\theta
_y + y'$, with $\Phi^T y'=0$. Then
%
\begin{eqnarray*}
\ud P_{\theta}(y) &=& c_1 \exp\biggl\{ -\frac{1}{2\sigma_n^2} (\|
\Phi\theta\|^2 + \|\Phi\theta_y\|^2 + \|y'\|^2 - 2 \theta^T \Phi
^T \Phi
\theta_y ) \biggr\}, \\[-2pt]
\ud\widetilde{W}(\theta) &=& c_2 \exp\biggl\{ -\frac{1}{2\tau_n^2} \|
\Phi\theta\|^2 \biggr\}, \\[-2pt]
\ud P_{\theta}(y) \,\ud\widetilde{W}(\theta) &=& c_1 c_2 \exp\biggl\{
-\frac{\sigma_n^2+\tau_n^2}{2\sigma_n^2 \tau_n^2} \biggl\|\Phi
\biggl(\theta- \frac{\tau_n^2}{\sigma_n^2+\tau_n^2} \theta_y\biggr)
\biggr\|^2 \\[-2pt]
&&\hphantom{c_1 c_2 \exp\biggl\{}
{} - \frac{1}{2(\sigma_n^2+\tau_n^2)} \|\Phi\theta
_y\|^2 -\frac{1}{2\sigma_n^2} \|y'\|^2 \biggr\},
\end{eqnarray*}
where $c_1 = (2\pi)^{-n/2} \sigma_n^{-n}$ and $c_2 = (2\pi)^{-k_n/2}
\tau_n^{-k_n} \det(\Phi^T \Phi)^{-1}$.

Using the Bayes rule, we get the density of $\widetilde{W}(\ud\theta|
Y)$, in which we recognize the normal distribution
%
%
\begin{equation} \label{eqexact-posterior}
\widetilde{W}(\ud\theta| Y) = \mathcal{N}\biggl(\frac{\tau_n^2}{\sigma
_n^2+\tau_n^2} \theta_Y, \frac{\sigma_n^2 \tau_n^2}{\sigma
_n^2+\tau
_n^2} (\Phi^T \Phi)^{-1}\biggr).
\end{equation}

So we have an exact expression for $\widetilde{W}(\ud\theta| Y)$, but
the centering and the variance do not correspond to the limit
distribution given in Theorem \ref{thmall-Gaussian-BvM}. Therefore, we
make use of the triangle inequality, with intermediate distribution
$Q=\mathcal{N}(\frac{\tau_n^2}{\sigma_n^2+\tau_n^2} \theta_Y,
\sigma_n^2 (\Phi^T \Phi)^{-1})$:
%
%
\begin{eqnarray} \label{eqgaussian-posterior}
&&
\bigl\| \widetilde{W}(\ud\theta| Y) - \mathcal{N}(\theta_Y,
\sigma_n^2 (\Phi^T \Phi)^{-1}) \bigr\|_{\mathrm{TV}} \nonumber\\[-9pt]\\[-9pt]
&&\qquad\leq
\bigl\| \widetilde{W}(\ud\theta| Y) - Q \bigr\|_{\mathrm{TV}} +
\| Q - \mathcal{N}(\theta_Y, \sigma_n^2 (\Phi^T \Phi)^{-1})
\|_{\mathrm{TV}}.
\nonumber
\end{eqnarray}

We first deal with the change in the variance, that is, the first term
on the right in (\ref{eqgaussian-posterior}).

Let $\alpha_n= \frac{\tau_n}{\sigma_n} \sqrt{\ln(1+\frac{\sigma
_n^2}{\tau_n^2})}$, and $f$ and $g$ be, respectively, the density
functions of $\mathcal{N}(0, I_{k_n})$\vadjust{\goodbreak} and $\mathcal{N}
(0, \frac{\tau_n^2}{\sigma_n^2+\tau_n^2} I_{k_n} )$. Let $U$ be a
random variable following the chi-square distribution with $k_n$
degrees of freedom $\chi^2(k_n)$.
Let $\sqrt{\Phi^T \Phi}$ be a square root of the matrix $\Phi^T
\Phi$.
The total variation norm is invariant under the bijective affine map
$\theta\mapsto\frac{1}{\sigma_n} \sqrt{\Phi^T \Phi} (\theta-
\frac{\tau_n^2}{\sigma_n^2+\tau_n^2} \theta_Y)$, so
%
\begin{eqnarray*}
\bigl\| \widetilde{W}(\ud\theta| Y) - Q \bigr\|_{\mathrm{TV}} &=&
\biggl\| \mathcal{N}(0, I_{k_n}) - \mathcal{N}\biggl(0, \frac
{\tau_n^2}{\sigma_n^2+\tau_n^2} I_{k_n}\biggr) \biggr\|_{\mathrm{TV}} \\[-2pt]
&=& \int_{\R^{k_n}} (g-f)_{+} =\int_{\|x\|\leq\sqrt
{k_n}\alpha_n} \bigl( g(x)-f(x) \bigr) \,\ud^n x \\[-2pt]
&=& \Proba\biggl(U\leq\frac{\sigma_n^2+\tau_n^2}{\tau_n^2} k_n\alpha
_n^2\biggr) - \Proba(U\leq k_n\alpha_n^2) \\[-2pt]
&=& \Proba\biggl(\sqrt{k_n} (\alpha_n^2 -1) \leq\frac
{U-k_n}{\sqrt{k_n}} \leq\sqrt{k_n} \biggl(\frac{\sigma_n^2+\tau
_n^2}{\tau_n^2} \alpha_n^2 -1\biggr) \biggr). 
\end{eqnarray*}
As $n$ goes to infinity, $\frac{U-k_n}{\sqrt{k_n}}$ 
converges toward $\mathcal{N}(0, 1 )$ in distribution. Using
the Taylor expansion of $\ln$, we find\vspace*{-3pt}
\[
\alpha_n^2 
= 1 - \frac{\sigma_n^2}{2\tau_n^2} + o\biggl(\frac{\sigma_n^2}{\tau
_n^2}\biggr)\vspace*{-3pt}
\]
and, therefore,
\begin{eqnarray*}
\sqrt{k_n} (\alpha_n^2 -1) &\sim& -\sqrt{k_n} \frac{\sigma
_n^2}{2\tau_n^2}, \\[-2pt]
\sqrt{k_n} \biggl(\frac{\sigma_n^2+\tau_n^2}{\tau_n^2} \alpha_n^2
-1\biggr) &\sim& \sqrt{k_n} \frac{\sigma_n^2}{2\tau_n^2}.
\end{eqnarray*}
Since $k_n = o(\tau_n^4/\sigma_n^4)$, both these quantities go to $0$.
As a consequence, $\| \widetilde{W}(\ud\theta|\break Y) - Q \|
_{\mathrm{TV}}$ goes to zero as $n$ goes to infinity.

Let us now deal with the centering term, that is, the second term on
the right in~(\ref{eqgaussian-posterior}).\vspace*{-3pt}
\begin{lem} \label{lemTV-gaussian-translation}
Let $U$ be a standard normal random variable, let $k\geq1$ and let
$Z\in\R^k$. Then
\[
\| \mathcal{N}(0, I_k) - \mathcal{N}(Z, I_k)
\|_{\mathrm{TV}} = \Proba( |U| \leq\|Z\|/2) \leq\|Z\|/\sqrt{2\pi
}.
\]
\end{lem}
\begin{pf}
Let $g$ be the density of $\mathcal{N}(0, I_k)$. Then
\begin{eqnarray*}
\| \mathcal{N}(0, I_k) - \mathcal{N}(Z, I_k
) \|_{\mathrm{TV}} &=& \int_{\R^k} \bigl(g(x)-g(x-Z)\bigr)_{+}
\,\ud^k x \\[-2pt]
&=& \int_{\{2 x^T Z \leq\|Z\|^2\}} \bigl(g(x)-g(x-Z)\bigr)\, \ud^k x
\\[-2pt]
&=& \Proba( U \leq\|Z\|/2) - \Proba( U + \|Z\| \leq\|Z\|/2) \\[-2pt]
&\leq& \|Z\|/\sqrt{2\pi}.
\end{eqnarray*}
The last line comes from the density of $\mathcal{N}(0, 1)$
being bounded by $1/\sqrt{2\pi}$.\vadjust{\goodbreak}~%
\end{pf}

Using again the invariance of the total variation norm under 
the bijective affine map $\theta\mapsto\frac{1}{\sigma_n} \sqrt
{\Phi^T
\Phi} (\theta- \frac{\tau_n^2}{\sigma_n^2+\tau_n^2} \theta_Y)$,
\begin{eqnarray*}
\| \mathcal{N}(\theta_Y, \sigma_n^2 (\Phi^T \Phi)^{-1}
) - Q \|_{\mathrm{TV}} &=& \biggl\| \mathcal{N}(0,
I_{k_n}) - \mathcal{N}\biggl(\frac{\sigma_n \sqrt{\Phi^T \Phi}
\theta_Y}{\tau_n^2+\sigma_n^2}, I_{k_n}\biggr) \biggr\|_{\mathrm{TV}} \\
&\leq& \frac{1}{\sqrt{2\pi}} \frac{\sigma_n}{(\tau_n^2+\sigma_n^2)}
\|\Phi\theta_Y\| \\
&\leq& \frac{1}{\sqrt{2\pi}} \frac{\sigma_n}{(\tau_n^2+\sigma_n^2)}
\bigl(\|F_0\| + \sqrt{\varepsilon^T \covmat\varepsilon}\bigr).
\end{eqnarray*}
$\varepsilon^T \covmat\varepsilon$ is a random variable following
$\sigma_n^2\chi^2(k_n)$ distribution. By Jensen's inequality, $\E
[\sqrt{\varepsilon^T \covmat\varepsilon}] \leq\sqrt{\E
[\varepsilon^T \covmat\varepsilon]} = \sigma_n\sqrt{k_n}$. Therefore,
\[
\E\| \mathcal{N}(\theta_Y, \sigma_n^2 (\Phi^T \Phi
)^{-1}) - Q \|_{\mathrm{TV}} \leq\frac{1}{\sqrt{2\pi}}
\frac{\sigma_n}{\tau_n^2+\sigma_n^2} \bigl(\|F_0\| + \sigma_n\sqrt
{k_n}\bigr),
\]
which goes to zero under the assumptions of Theorem \ref{thmall-Gaussian-BvM}.

To conclude the proof, note that we deduce the results on $W(\ud F |
Y)$ from the ones on $\widetilde{W}(\ud\theta| Y)$, by the linear
relation $F=\Phi\theta$.

\subsection{\texorpdfstring{Proof of Theorem \protect\ref{thmsmooth-gaussian-BvM}}{Proof of Theorem 2}}
\label
{secproof-smooth-gaussian}
We make the proof for $\widetilde{W}(\ud\theta| Y)$. Then the result
for $W(\ud F | Y)$ is immediate. Our method is adapted from \citet
{BoucheronGassiat09}.

To any probability measure $P$ on $\R^{k_n}$, we associate the probability
%
%
\begin{equation} \label{eqdef-P-M}
P^M 
= \frac{P(\cdot\cap\mathcal{E}_{\theta_0, \Phi}(M)
)}{P( \mathcal{E}_{\theta_0, \Phi}(M) )}
\end{equation}
with support in $\mathcal{E}_{\theta_0, \Phi}(M)$. It can be easily
checked that
%
%
\begin{equation} \label{eqtruncated-distribution}
\| P - P^M \|_{\mathrm{TV}} = P(\mathcal{E}_{\theta_0,
\Phi}^c (M)).
\end{equation}
%

\begin{sloppypar}
The proof is divided into three steps based on the use of $M_n$ as a
threshold to truncate the probability distributions. Lemma \ref
{lemgaussian-part} below controls $\E\| \mathcal{N}(\theta
_Y, \sigma_n^2 (\Phi^T \Phi)^{-1}) - \mathcal{N}^{M_n}(\theta
_Y, \sigma_n^2 (\Phi^T \Phi)^{-1}) \|_{\mathrm{TV}}$,
Lemma \ref{lemcentral-part} controls $\E\| \widetilde
{W}^{M_n}(\ud\theta| Y) - \mathcal{N}^{M_n}(\theta_Y, \sigma_n^2
(\Phi^T \Phi)^{-1}) \|_{\mathrm{TV}}$ and Proposition \ref
{propposterior-part} controls $\E\| \widetilde{W}(\ud\theta|
Y) - \widetilde{W}^{M_n}(\ud\theta| Y) \|_{\mathrm{TV}}$. Taken
together, these results give Theorem~\ref{thmsmooth-gaussian-BvM}.
\end{sloppypar}
\begin{lem} \label{lemgaussian-part} 
If $k_n<4 M_n$, then
\[
\E\| \mathcal{N}(\theta_Y, \sigma_n^2 (\Phi^T \Phi
)^{-1}) - \mathcal{N}^{M_n}(\theta_Y, \sigma_n^2 (\Phi^T \Phi
)^{-1}) \|_{\mathrm{TV}} \leq2 e^{-{(\sqrt
{M_n}-2\sqrt{k_n})^2}/{8} }.
\]
\end{lem}

If $k_n=o(M_n)$, for $n$ large enough, this bound can be replaced
by\vadjust{\goodbreak}
$e^{-M_n/9}$.
\begin{pf*}{Proof of Lemma \ref{lemgaussian-part}}
To control this quantity, we consider two cases, 
depending on whether $\theta_Y$ is near or far from $\theta_0$:
%
%
\begin{eqnarray}\label{eqgaussian-truncation}
&&\| \mathcal{N}(\theta_Y, \sigma_n^2 (\Phi^T \Phi)^{-1}
) - \mathcal{N}^{M_n}(\theta_Y, \sigma_n^2 (\Phi^T \Phi)^{-1}
) \|_{\mathrm{TV}} \nonumber\\[2pt]
&&\qquad = \mathcal{N}(\theta_Y, \sigma_n^2 (\Phi^T \Phi
)^{-1}) (\mathcal{E}_{\theta_0, \Phi}^c (M_n) ) \nonumber\\[-7pt]\\[-7pt]
&&\qquad\leq\mathbh{1}_{(\theta_Y-\theta_0)^T \Phi^T\Phi(\theta
_Y-\theta_0) > \sigma_n^2 M_n/4} \nonumber\\[2pt]
&&\qquad\quad{} + \mathcal{N}(\theta_0, \sigma_n^2 (\Phi^T \Phi)^{-1}
) \bigl(\mathcal{E}_{\theta_0, \Phi}^c (M_n/4) \bigr).
\nonumber
\end{eqnarray}
Let $U$ be a random variable following a $\chi^2(k_n)$ distribution.
Taking the expectation on both sides of (\ref{eqgaussian-truncation}) gives
\[
\E\| \mathcal{N}(\theta_Y, \sigma_n^2 (\Phi^T \Phi
)^{-1}) - \mathcal{N}^{M_n}(\theta_Y, \sigma_n^2 (\Phi^T \Phi
)^{-1}) \|_{\mathrm{TV}} \leq2 \Proba(U>M_n/4).
\]
Now, Cirelson's inequality [see, e.g., \citet{Massart2007}]
%
%
\begin{equation} \label{eqCirelson}
\Proba\bigl(\sqrt{U}>\sqrt{k_n}+\sqrt{2 x}\bigr) \leq\exp(-x)
\end{equation}
used with $x=\frac{(\sqrt{M_n}-2\sqrt{k_n})^2}{8}$ implies
Lemma \ref{lemgaussian-part}.
\end{pf*}
\begin{lem} \label{lemcentral-part}
If $ \sup_{\|\Phi h\|^2 \leq\sigma_n^2 M_n, \|\Phi g\|^2
\leq\sigma_n^2 M_n} \frac{w(\theta_0+h)}{w(\theta_0+g)}
\rightarrow1$
as $n\rightarrow\infty$, then
\[
\E\bigl\| \widetilde{W}^{M_n}(\ud\theta| Y) - \mathcal{N}^{M_n}
(\theta_Y, \sigma_n^2 (\Phi^T \Phi)^{-1}) \bigr\|_{\mathrm{TV}}
\rightarrow0 \qquad\mbox{as } n\rightarrow\infty.
\]
\end{lem}
\begin{pf}
Let us first note that, for every $\theta$ and $\tau$ in $\R^{k_n}$,
for every $Y\in\R^n$,
%
%
\begin{eqnarray} \label{eqexact-LAN}
\frac{\ud P_\theta(Y)}{\ud P_\tau(Y)} &=& \exp\biggl\{ \frac{-\|\Phi
\theta\|^2+\|\Phi\tau\|^2-2 Y^T \Phi(\tau-\theta)}{2\sigma_n^2}
\biggr\} \nonumber\\[-7pt]\\[-7pt]
&=& \frac{\ud\mathcal{N}(\theta_Y, \sigma_n^2 (\Phi^T \Phi
)^{-1})(\theta)}{\ud\mathcal{N}(\theta_Y, \sigma_n^2 (\Phi^T
\Phi)^{-1})(\tau)}.\nonumber
\end{eqnarray}
This directly comes from the expressions for the Gaussian densities.

In the following the first lines are just rewriting. Then we use
Jensen's inequality with the convex function\vadjust{\goodbreak} $x
\mapsto(1-x)_{+}$, and make use of (\ref{eqexact-LAN}). We abbreviate
$\mathcal{N}^{M_n} (\theta_Y, \sigma_n^2 (\Phi^T \Phi)^{-1})$ into
$\mathcal{N}^{M_n}$:
\begin{eqnarray*}
&&\bigl\| \widetilde{W}^{M_n}(\ud\theta| Y) - \mathcal{N}^{M_n}
\bigr\|_{\mathrm{TV}} \\[2pt] 
%
&&\qquad= \int\biggl(1- \frac{ \ud\mathcal{N}^{M_n}(\theta) }{ \ud\widetilde
{W}^{M_n}(\theta| Y) } \biggr)_{+} \,\ud\widetilde{W}^{M_n}(\theta|
Y) \\[2pt]
&&\qquad= \int\biggl(1- \frac{ \ud\mathcal{N}^{M_n}(\theta) \int
({w(\tau
)}/{\ud\mathcal{N}^{M_n}(\tau)}) \,\ud P_\tau(Y) \,\ud\mathcal
{N}^{M_n}(\tau) }{ w(\theta) \,\ud P_\theta(Y) } \biggr)_{+} \,\ud
\widetilde{W}^{M_n}(\theta| Y) \\[2pt]
&&\qquad\leq\int\!\int\biggl(1-\frac{w(\tau) \,\ud\mathcal{N}^{M_n}(\theta)
\,\ud P_\tau(Y)}{w(\theta) \,\ud\mathcal{N}^{M_n}(\tau) \,\ud P_\theta
(Y)}\biggr)_{+} \,\ud\mathcal{N}^{M_n}(\tau) \,\ud\widetilde
{W}^{M_n}(\theta| Y) \\[2pt]
&&\qquad= \int\!\int\biggl(1-\frac{w(\tau)}{w(\theta)}\biggr)_{+}
\,\ud\mathcal
{N}^{M_n}(\tau) \,\ud\widetilde{W}^{M_n}(\theta| Y) \\[2pt] 
&&\qquad\leq1- \inf_{\|\Phi h\|^2 \leq\sigma_n^2 M_n, \|\Phi g\|^2 \leq
\sigma_n^2 M_n} \frac{w(\theta_0+h)}{w(\theta_0+g)}.
\end{eqnarray*}
\upqed\end{pf}
\begin{prop}[(Posterior concentration)] \label{propposterior-part}
Suppose that conditions \ref{hypcontinuous-w}, \ref
{hypk-n-M-n} and \ref{hypweight} of Theorem \ref
{thmsmooth-gaussian-BvM} hold. Then
\begin{eqnarray*}
\E\bigl\| \widetilde{W}(\ud\theta| Y) - \widetilde{W}^{M_n}(\ud
\theta| Y) \bigr\|_{\mathrm{TV}} &=& \E[ \widetilde{W}(
\mathcal{E}_{\theta_0, \Phi}^C (M_n) | Y )] \\[-1pt]
&\rightarrow& 0 \qquad\mbox{as } n\rightarrow\infty.
\end{eqnarray*}
\end{prop}

Proposition \ref{propposterior-part} is proved in Appendix A %
in the supplemental article [\citet{BvMGRegSupplement}]. However, we state
here the following important lemma, because of its significance.
\begin{lem} \label{leminvariant-translation}
Let $a\in\R^n$ such that $\Phi^T a=0$. Then, for any $y\in\R^n$,
$W(\ud F |\allowbreak Y=y) = W(\ud F |Y=y+a)$.
\end{lem}

Lemma \ref{leminvariant-translation} states that the distribution
$W(\ud F |Y)$ is invariant under any translation of $Y$ orthogonal to
$\phispan$. Now, regard $W(\ud F |Y)$ as a random variable. Then any
statement on $W(\ud F |Y)$ or $\widetilde{W}(\ud\theta| Y)$ valid
when $Y\sim\mathcal{N}(F_0, \sigma_n^2 I_n)$ with $F_0 \in\phispan$
can be extended at zero cost by Lemma \ref{leminvariant-translation}
to the case $F_0\in\R^n$. For instance, proving Proposition \ref
{propposterior-part} in the case $F_0=\Phi\theta_0$ is enough.\vspace*{3pt}

\subsection{\texorpdfstring{Proof of Theorem \protect\ref{thmfunctional-base}}{Proof of Theorem 3}}
\label{secproof-functional-all-gaussian}

We begin with (\ref{eqbayesian-functionals}). Consider the following
Taylor expansion:
\begin{eqnarray*}
&&G(F)-G\bigl(Y_\phispan\bigr)\\[-1pt]
&&\qquad= \dot{G}_{F_\phispan} \bigl(F-Y_\phispan\bigr) \\[-1pt]
&&\qquad\quad{} + \frac{1}{2} \int_0^1 (1-t) D_{F_\phispan+t (F-F_\phispan)}^2
G\bigl(F-F_\phispan, F-F_\phispan\bigr) \,\ud t \\[-1pt]
&&\qquad\quad{} - \frac{1}{2} \int_0^1 (1-t) D_{F_\phispan+t (Y_\phispan
-F_\phispan)}^2 G\bigl(Y_\phispan-F_\phispan, Y_\phispan-F_\phispan\bigr)
\,\ud t
\end{eqnarray*}
using the Lagrange form of the error term. Suppose that $F\in\phispan$,
$\|F-F_\phispan\|^2\leq\sigma_n^2 M_n$ and $\|Y_\phispan
-F_\phispan\|
^2\leq\sigma_n^2 M_n$. Then, for any $b\in\R^p$,
\[
\bigl| b^T \bigl(G(F)-G\bigl(Y_\phispan\bigr) - \dot{G}_{F_\phispan}
\bigl(F-Y_\phispan\bigr) \bigr) \bigr| \leq\|b\| B_{F_\phispan}(M_n).
\]
On the other hand, $\sqrt{b^T\Gamma_{F_\phispan} b} \geq\sqrt{\|
\Gamma_{F_\phispan}^{-1}\|^{-1}} \|b\|$. Moreover,
%
\begin{eqnarray*}
&&\biggl\| W\biggl(\ud\frac{b^T \dot{G}_{F_\phispan} (F-Y_\phispan
)}{\sqrt{b^T \Gamma_{F_\phispan} b}} \Big| Y\biggr) - \mathcal{N}
(0, 1) \biggr\|_{\mathrm{TV}} \\
&&\qquad\leq\bigl\| W(\ud F | Y) - \mathcal{N}\bigl(Y_\phispan, \sigma_n^2
\covmat\bigr) \bigr\|_{\mathrm{TV}}.
\end{eqnarray*}

Let $\eta_n = \sqrt{\|\Gamma_{F_\phispan}^{-1}\|}
B_{F_\phispan}(M_n)$, which tends to $0$ by hypothesis. Let also
\[
I_{\eta_n} = \{ x\in\R\dvtx\exists x' \in I, |x-x'|\leq\eta_n \}.
\]
Note that $\psi(I_{\eta_n}) \leq\psi(I) + \sqrt{\frac{2}{\pi
}} \eta_n$.

Gathering all this information, we can get the upper bound
\begin{eqnarray*}
&&W\biggl( \frac{b^T (G(F)-G(Y_\phispan))}{\sqrt
{b^T\Gamma_{F_\phispan} b}} \in I \Big| Y\biggr) \\
&&\qquad\leq W\biggl( \frac{b^T \dot{G}_{F_\phispan} (F-Y_\phispan
)}{\sqrt{b^T \Gamma_{F_\phispan} b}} \in I_{\eta_n} \Big| Y\biggr) \\
&&\qquad\quad{} + \mathbh{1}_{\|Y_\phispan-F_\phispan\|^2> \sigma_n^2 M_n} + W
\bigl(\bigl\|F-F_\phispan\bigr\|^2> \sigma_n^2 M_n | Y\bigr) \\
&&\qquad\leq\psi(I) + \sqrt{\frac{2}{\pi}} \eta_n + \bigl\| W(\ud F | Y) -
\mathcal{N}\bigl(Y_\phispan, \sigma_n^2 \covmat\bigr) \bigr\|_{\mathrm{TV}} \\
&&\qquad\quad{} + \mathbh{1}_{\|Y_\phispan-F_\phispan\|^2> \sigma_n^2 M_n} + W
\bigl(\bigl\|F-F_\phispan\bigr\|^2> \sigma_n^2 M_n | Y\bigr).
\end{eqnarray*}

A lower bound is obtained in the same way. Taking the expectation,
%
%
\begin{eqnarray} \label{eqfunctionnal-deviation}
&&\E\biggl| W\biggl( \frac{b^T (G(F)-G(Y_\phispan)
)}{\sqrt{b^T\Gamma_{F_\phispan} b}} \in I\Big | Y\biggr) - \psi(I)
\biggr| \nonumber\\
&&\qquad\leq o(1) + \Proba\bigl(\bigl\|Y_\phispan-F_\phispan\bigr\|^2> \sigma_n^2 M_n\bigr) \\
&&\qquad\quad{} + \E\bigl[W\bigl(\bigl\|F-F_\phispan\bigr\|^2> \sigma_n^2 M_n
| Y\bigr)\bigr].
\nonumber
\end{eqnarray}
But $\|Y_\phispan-F_\phispan\|^2$ follows the $\sigma_n^2\chi^2(k_n)$
distribution, and since $k_n = o(M_n)$,
\[
\Proba\bigl(\bigl\|Y_\phispan-F_\phispan\bigr\|^2> \sigma_n^2 M_n\bigr) = o(1).
\]

To bound (\ref{eqfunctionnal-deviation}), we use the following:
\begin{lem} \label{lemposterior-consistency-all-gaussian}
Suppose that the conditions of either Theorems \ref
{thmall-Gaussian-BvM} or \ref{thmsmooth-gaussian-BvM} are
satisfied. Then
\[
\E\bigl[W\bigl(\bigl\|F-F_\phispan\bigr\|^2> \sigma_n^2 M_n | Y
\bigr)\bigr] \rightarrow0 \qquad\mbox{as } n \rightarrow\infty.
\]
\end{lem}
\begin{pf}
For smooth priors, this is an immediate corollary of Proposition~\ref
{propposterior-part}. Let us suppose we are under the conditions of
Theorem \ref{thmall-Gaussian-BvM}.

Let $Z$ be a $\mathcal{N}(0, \frac{\sigma_n^2 \tau_n^2}{\sigma
_n^2+\tau_n^2} \covmat)$ random vector in $\R^n$ independent on
$Y$, and $U$ a random variable following $\chi^2(k_n)$. From (\ref
{eqexact-posterior}) we get
\[
W(\ud F | Y) = \mathcal{N}\biggl(\frac{\tau_n^2}{\sigma_n^2+\tau_n^2}
Y_\phispan, \frac{\sigma_n^2 \tau_n^2}{\sigma_n^2+\tau_n^2}
\covmat
\biggr).
\]
Therefore,
\begin{eqnarray*}
&&W\bigl(\bigl\|F-F_\phispan\bigr\|^2> \sigma_n^2 M_n | Y\bigr) \\
&&\qquad = \Proba\biggl( \biggl\| Z + \frac{\tau_n^2}{\sigma_n^2+\tau_n^2}
Y_\phispan- F_\phispan\biggr\|^2 > \sigma_n^2 M_n\biggr) \\
&&\qquad \leq\Proba\biggl( \| Z \| > \sigma_n \sqrt{M_n} - \biggl\|
\frac{\tau_n^2}{\sigma_n^2+\tau_n^2} Y_\phispan- F_\phispan\biggr\|
\biggr) \\
&&\qquad\leq\cases{
1, \qquad \mbox{if $\displaystyle \biggl\|\frac{\tau_n^2}{\sigma_n^2+\tau_n^2}
Y_\phispan- F_\phispan\biggr\| > \frac{2 \sigma_n \sqrt{M_n}}{3}$}, \vspace*{2pt}\cr
\displaystyle \Proba\biggl(\|Z\|^2 > \sigma_n^2 \frac{M_n}{9} \biggr) = \Proba\biggl(U
> \frac{\sigma_n^2+\tau_n^2}{\tau_n^2} \frac{M_n}{9} \biggr), \cr
\hspace*{33.8pt}\mbox{otherwise.}}
\end{eqnarray*}
Since $k_n=o(M_n)$, $\Proba(U>M_n/9) = o(1)$. On the other hand,
\begin{eqnarray*}
\biggl\|\frac{\tau_n^2}{\sigma_n^2+\tau_n^2} Y_\phispan- F_\phispan
\biggr\| &=& \biggl\|\covmat\biggl(\frac{\tau_n^2}{\sigma_n^2+\tau_n^2}
\varepsilon+ \frac{\sigma_n^2}{\sigma_n^2+\tau_n^2} F_0\biggr)
\biggr\| \\
&\leq&\|\covmat\varepsilon\| + \frac{\sigma_n}{\sqrt{\sigma
_n^2+\tau
_n^2}} \|F_0\|.
\end{eqnarray*}
Since $\|F_0\| = o(\tau_n^2/\sigma_n)$, $\frac{\sigma_n^2 \|F_0\|
^2}{\sigma_n^2+\tau_n^2} = o(1) < \frac{M_n}{9}$ for $n$ large enough.
$\|\covmat\varepsilon\|^2$ is a~$\sigma_n^2 \chi^2(k_n)$ variable.
Therefore, for $n$ large enough,
\[
\E\bigl[W\bigl(\bigl\|F-F_\phispan\bigr\|^2> \sigma_n^2 M_n | Y
\bigr)\bigr] \leq2 \Proba(U>M_n/9) = o(1).
\]
\upqed\end{pf}
Now, (\ref{eqfunctionnal-deviation}) gives
(\ref{eqbayesian-functionals}).\vadjust{\goodbreak}

The proof of the frequentist assertion (\ref
{eqfrequentist-functionals}) is similar and delayed to Appendix C %
in the supplemental article [\citet{BvMGRegSupplement}].

\section*{Acknowledgments}

The author would like to thank E. Gassiat and I.~Cas\-tillo for valuable
discussions and suggestions.

\begin{supplement}[id=suppA]
\stitle{Supplement to ``Bernstein--von Mises theorems for Gaussian regression with
increasing number of regressors''\vadjust{\goodbreak}}
\slink[doi]{10.1214/11-AOS912SUPP} 
\sdatatype{.pdf}
\sfilename{aos912\_supp.pdf}
\sdescription{This
contains the proofs of various technical results stated in the main
article ``Bernstein--von Mises Theorems for Gaussian regression with
increasing number of regressors.''}
\end{supplement}

%

\printaddresses

\end{document}